\documentclass[a4paper]{scrartcl}

\usepackage[english]{babel}
\usepackage[utf8x]{inputenc}
\usepackage[T1]{fontenc}
\usepackage{amsmath}

\usepackage{tkz-euclide}
\usepackage{amsfonts}
\usepackage{enumerate}
\usepackage{float}
\usetikzlibrary{calc}
\usepackage{MnSymbol,blindtext,wasysym,etoolbox}
\usepackage{pgfplots}
\usepackage{multicol}
\usepackage{comment}

\usepackage[a4paper,top=3cm,bottom=2cm,left=3cm,right=3cm,marginparwidth=1.75cm]{geometry}

\usepackage{amsmath}
\usepackage{amsmath, dsfont,mathtools,booktabs}
\newcommand*\diff{\mathop{}\!\mathrm{d}}

\usepackage{mathtools}
\usepackage{enumitem}
\usepackage{bbm}
\usepackage{titlesec}
\usepackage{caption}
\usepackage{subcaption}
\usepackage{verbatim}
\usepackage{graphicx}
\usepackage[colorinlistoftodos]{todonotes}
\usepackage[colorlinks=true, allcolors=blue]{hyperref}
\usepackage{algorithm}
\usepackage[noend]{algpseudocode}
\usepackage{mathtools}
\usepackage[thinlines]{easytable}
\usepackage{authblk}
\allowdisplaybreaks

\makeatletter
\def\BState{\State\hskip-\ALG@thistlm}
\makeatother
\usetikzlibrary{snakes}

\newtheorem{definition}{Definition}
\newtheorem{problem}{Problem}

\newtheorem{remark}{Remark}
\newtheorem{configuration}{Configuration}

\title{\textbf{Adaptive time-step control for a monolithic multirate scheme coupling the heat and wave equation}}
\author[1]{Martyna Soszy\'nska\thanks{martyna.soszynska@ovgu.de}}
\author[1]{Thomas Richter\thanks{thomas.richter@ovgu.de}}
\affil[1]{Institut f\"ur Analysis und Numerik, Otto-von-Guericke Univerist\"at}
\date{\today}
\begin{document}
\maketitle
\section*{Abstract}
\label{abstract}

  We consider the dynamics of a parabolic and a hyperbolic equation coupled on a common interface and develop time-stepping schemes that can use different time-step sizes for each of the subproblems. The problem is formulated in a strongly coupled (monolithic) space-time framework. Coupling two different step sizes monolithically gives rise to large algebraic systems of equations where multiple states of the subproblems must be solved at once. For efficiently solving these algebraic systems, we inherit ideas from the partitioned regime and present two decoupling methods, namely a partitioned relaxation scheme and a shooting method.
  
  Furthermore, we develop an a posteriori error estimator serving as a mean for an adaptive time-stepping procedure. The goal is to optimally balance the time step sizes of the two subproblems. The error estimator is based on the dual weighted residual method and relies on the space-time Galerkin formulation of the coupled problem.

  As an example, we take a linear set-up with the heat equation coupled to the wave equation. We formulate the problem in a monolithic manner using the space-time framework. In numerical test cases, we demonstrate the efficiency of the solution process and we also validate the accuracy of the a posteriori error estimator and its use for controlling the time step sizes.

\section{Introduction}
\label{introduction}

In this work, we are going to work with surface coupled multiphysics problems that are inspired by fluid-structure interaction (FSI) problems~\cite{Richter2017}. We couple the heat equation with the wave equation through an interface, where the typical FSI coupling conditions or Dirichlet-Neumann type act. Despite of its simplicity, each of the subproblems exhibits different temporal dynamics which is also found in FSI. The solution of the heat equation, as a parabolic problem, manifests smoothing properties, thus it can be characterized as a problem with slow temporal dynamics. The wave equation, on the other hand, is an example of a hyperbolic equation with highly oscillatory properties. 

FSI problems are characterized by two specific difficulties: the coupling of an equation of parabolic type with one of hyperbolic type gives rise to regularity problems at the interface. Further, the added mass effect~\cite{CausinGerbeauNobile2005}, which is present for problems coupling materials of a similar density, calls for discretization and solution schemes which are strongly coupled. This is the monolithic approach for modeling FSI, in contrast to partitioned approaches, where each of the subproblems is treated and solved as a separate system. While the monolithic approach allows for a more rigorous mathematical setting and the use of large time steps, the partitioned approach allows using fully optimized separate techniques for both of the subproblems. Most realizations for FSI, such as the technique described here, have to be regarded as a blend of both philosophies: while the formulation and discretization are monolithic, ideas of partitioned approaches are borrowed for solving the algebraic problems. 

Featuring distinct time scales in each of the problems, the use of multirate time-stepping schemes with adapted step sizes for fluid and solid is obvious. For parabolic problems, the concept of multirate time-stepping was discussed in \cite{Dawson1991}, \cite{Blum1992} and \cite{Faille2009}. In the hyperbolic setting, it was considered in~\cite{BergerMarsha1985}, \cite{Collino2003part1} \cite{Collino2003part2} and \cite{Piperno2006}. In the context of fluid-structure interactions, such subcycling methods are used in aeroelasticity~\cite{Piperno1997}, where explicit time integration schemes are used for the flow problem and implicit schemes for the solid problem~\cite{DeMoerlooseetal2018}. In the low Reynolds number regime, common in hemodynamics, the situation is different. Here, implicit and strongly coupled schemes are required by the added mass effect. Hence, large time steps can be applied for the flow problem, but smaller time steps might be required within the solid. A study on benchmark problems in fluid dynamics (Sch\"afer, Turek '96~\cite{SchaeferTurek1996}) and FSI presented in~\cite{HronTurek2006} shows that FSI problems demand a much smaller step size, although the problem configuration and the resulting nonstationary dynamics are very similar to oscillating solutions with nearly the same period~\cite{RichterWick2015_time}.

We will derive a monolithic variational formulation for FSI like problems that can handle different time step sizes in the two subproblems. Implicit coupling of two problems with different step sizes will give rise to very large systems where multiple states must be solved at once.  In Section~\ref{decoupling_methods} we will study different approaches for an efficient solution of these coupled systems, a simple partitioned relaxation scheme and a shooting like approach.

Next, in Section~\ref{goal_oriented_estimation} we present a posteriori error estimators based on the dual weighted residual method~\cite{BeckerRannacher2001} for automatically identifying optimal step sizes for the two subproblems. Numerical studies on the efficiency of the time adaptation procedure are presented in Section \ref{numerical_results}.

%
\section{Presentation of the model problem}
\label{presentation_of_the_problem}
Let us consider the time interval $I = [0, T]$ and two rectangular domains $\Omega_f = (0, 4) \times (0, 1)$, $\Omega_s = (0, 4) \times (0, -1)$. The interface is defined as $\Gamma \coloneqq \overline{\Omega}_f \cap \overline{\Omega}_s =  (0,4)\times \{0\}$. The remaining boundaries are determined as $\Gamma_f^1 \coloneqq \{0 \} \times (0, 1)$, $\Gamma_f^2 \coloneqq (0, 4) \times \{ 1\}$, $\Gamma_f^3 \coloneqq \{4 \} \times (0, 1)$ and $\Gamma_s^1 \coloneqq \{0 \} \times (-1, 0)$, $\Gamma_s^2 \coloneqq (0, 4) \times \{ -1\}$, $\Gamma_s^3 \coloneqq \{4 \} \times (-1, 0)$. The domain is illustrated in Figure \ref{domain}.
In the domain $\Omega_f$ we pose the heat equation parameterized by
the diffusion parameter $\nu > 0$ with an additional transport term
controlled by $\beta \in \mathds{R}^2$. In the domain $\Omega_s$ we
set the wave equation. By $\sqrt{\lambda}$ we denote the propagation
speed and by $\delta \geq 0$ a damping parameter. On the interface, we
set both kinematic and dynamic coupling conditions. The former
guarantees the continuity of displacement and velocity along the
interface. The latter establishes the balance of normal stresses. The
exact values of the parameters read as
\[
\nu = 0.001,\quad
\beta =  \left(\begin{matrix} 2 \\ 0 \end{matrix}\right),\quad
\lambda = 1000,\quad \delta = 0.1
\]
and the complete set of equations is given by
\begin{equation*}
\begin{cases}
  \partial_t v_f - \nu \Delta v_f + \beta \cdot \nabla v_f = g_f,\quad
- \Delta u_f = 0 & \textnormal{in } I \times \Omega_f, \\
\partial_t v_s - \lambda \Delta u_s - \delta \Delta v_s = g_s, \quad
\partial_t u_s = v_s & \textnormal{in } I \times \Omega_s, \\
u_f = u_s,\quad
v_f = v_s,\quad
\lambda \partial_{\vec{n}_s} u_s = -\nu \partial_{\vec{n}_f}v_f & \textnormal{on } I \times \Gamma, \\
u_f = v_f = 0 & \textnormal{on } I \times \Gamma_f^2,  \\
u_s = v_s = 0 & \textnormal{on } I \times \Gamma_s^1\cup \Gamma_s^3,  \\
u_f(0) = v_f(0) = 0 & \textnormal{in } \Omega_f, \\
u_s(0) = v_s(0) = 0 & \textnormal{in } \Omega_s \\
\end{cases}
\end{equation*}
We use symbols $\vec{n}_f$ and $\vec{n}_s$ to distinguish between
normal vectors for different space domains. 

\begin{figure}[t]
\begin{center}
\begin{tikzpicture}[scale = 0.85]
    \draw (-0.5, 0) -- (5.5, 0);
    \draw[ultra thick] (0,0.5) -- (0.5, 0.5);
    \draw[ultra thick] (0.5,0) -- (5, 0);
    \draw[fill = black] (0, 0.5) circle (0.075cm);
    \draw[fill = white] (0.5, 0.5) circle (0.075cm);
    \draw[fill = black] (0.5, 0) circle (0.075cm);
    \draw[fill = white] (5, 0) circle (0.075cm);
    \node at (0, -0.3) {0};
    \node at (0.5, -0.3) {0.1};
    \node at (5, -0.3) {1};
    \draw[white] (0.0, 1.0) -- (1.0, 1.0);
\end{tikzpicture}
\caption{Function $g_2$ on $I=[0,T)$ for $T = 1$. }
\label{g_2}
\label{proba}
\end{center}
\end{figure}
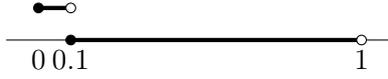

The external
forces are set to be products of functions of space and time
${g_f(\vec{x}, t) \coloneqq g_f^1(\vec{x})g^2(t)}$ and ${g_s(\vec{x}, t) \coloneqq
  g_s^1(\vec{x})g^2(t)}$ where $g_f^1$, $g_s^1$ are space components and $g^2$ is a time component. We will consider two configurations of the right hand side. In Configuration \ref{configuration_1}, the right hand side is concentrated in $\Omega_f$ where the space component consists of an exponential function centered around $\left(\frac{1}{2}, \frac{1}{2} \right)$. For Configuration \ref{configuration_2} we take a space component concentrated in $\Omega_s$ with an exponential function centered around $\left(\frac{1}{2}, -\frac{1}{2} \right)$.
\begin{configuration}\label{configuration_1}
  \begin{alignat*}{2}
g_f^1(\vec{x})\coloneqq &e^{-\left((x_1 - \frac{1}{2})^2 + (x_2 -
  \frac{1}{2})^2\right)}, \quad && \vec{x} \in \Omega_f \\ 
g_s^1(\vec{x})\coloneqq &0, \quad && \vec{x} \in \Omega_s
\end{alignat*}
\end{configuration}
\begin{configuration}\label{configuration_2}
\begin{alignat*}{2}
g_f^1(\vec{x})\coloneqq &0,\quad  && \vec{x} \in \Omega_f \\
g_s^1(\vec{x})\coloneqq &e^{-\left((x_1 - \frac{1}{2})^2 + (x_2 + \frac{1}{2})^2\right)}, \quad  && \vec{x} \in \Omega_s
\end{alignat*}
\end{configuration}
For both cases, we chose the same time component $g^2(t) \coloneqq \mathbbm{1}_{\big[\lfloor t \rfloor, \lfloor t \rfloor + \frac{1}{10}\big)}(t)$ for $t \in I$ illustrated in Figure \ref{g_2}.

\begin{figure}[t]
\begin{center}
\begin{tikzpicture}[scale = 1.0]
    \draw (0.0, 1.0) -- (4.0, 1.0) -- (4.0, -1.0);
    \draw (4.0, -1.0) -- (0.0, -1.0) -- (0.0, 1.0);
    \draw[dashed] (0.0, 0.0) -- (4.0, 0.0);
    \node at (0.35, 0.25) {$\Omega_f$};
    \node at (-0.35, 0.5) {$\Gamma_f^1$};
    \node at (-0.35, -0.5) {$\Gamma_s^1$};
    \node at (4.35, 0.5) {$\Gamma_f^3$};
    \node at (4.35, -0.5) {$\Gamma_s^3$};
    \node at (0.35, -0.75) {$\Omega_s$};
    \node at (2.0, 1.35) {$\Gamma_f^2$};
    \node at (2.0, 0.35) {$\Gamma$};
    \node at (2.0, - 1.35) {$\Gamma_s^2$};
\end{tikzpicture}
\caption{View of the domain split into ``fluid'' $\Omega_f$ and
  ``solid'' $\Omega_s$ along the common interface~$\Gamma$. }
\label{domain}
\end{center}
\end{figure}
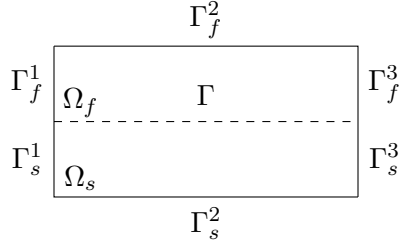

Since our example might be treated as a simplified case of an FSI problem, in the text we will use the corresponding nomenclature. We will refer to domain $\Omega_f$ as the \textit{fluid domain} and the problem defined there as the \textit{fluid problem}. Similarly, we will use \textit{solid domain} and \textit{solid problem} phrases.
\subsection{Continuous variational formulation}
As the first step, let us introduce a family of Hilbert spaces, which will be later on used as the trial and test spaces for our variational problems
\begin{equation*}
    X(V) = \left\{ v \in L^2(I, V) |\; \partial_t v \in L^2(I, V^*) \right\}.
\end{equation*}
Because we would like to incorporate the Dirichlet boundary conditions on $\Gamma_f^2$ and $\Gamma_s^1$, $\Gamma_s^3$ into spaces of solutions, for $\Upsilon \subset \partial \Omega$, we define
\begin{equation*}
    H^1_0(\Omega; \Upsilon) = \left\{v \in H^1(\Omega)|\; v_{|\Upsilon} = 0 \right\}.
\end{equation*}
Note that $\left(H^1_0(\Omega;\Upsilon)\right)^* = H^{-1}(\Omega)$.
For our example, we choose $H_f \coloneqq H^1_0(\Omega_f; \Gamma_f^2)$ and
$H_s \coloneqq H^1_0(\Omega_s;\Gamma_s^1\cup\Gamma_s^3)$ for representing space. We take $X_f \coloneqq (X(H_f))^2$, $X_s\coloneqq (X(H_s))^2$ and $X = X_f \times X_s$ for space-time trial and test function spaces. Below we present notations for inner products and duality pairings:
\begin{alignat*}{2}
& (u, \varphi)_f\coloneqq (u, \varphi)_{L^2(\Omega_f)}, \quad &&  \langle u, \varphi \rangle_f  \coloneqq \langle u, \varphi \rangle_{H^{-1}(\Omega_f) \times H_f}, \\
& (u, \varphi)_s\coloneqq (u, \varphi)_{L^2(\Omega_s)}, && \langle u, \varphi \rangle_s  \coloneqq \langle u, \varphi \rangle_{H^{-1}(\Omega_s) \times H_s}, \\
&  && \langle u, \varphi \rangle_{\Gamma} \coloneqq \langle u, \varphi \rangle_{H^{-\frac{1}{2}}(\Gamma) \times H^{\frac{1}{2}}(\Gamma)}
\end{alignat*}
To shorten the notation, we introduce the abbreviations
\[
\begin{aligned}
  \vec{U}_f &\coloneqq \left(\begin{matrix} 
    u_f \\ v_f \end{matrix}\right),& 
  \vec{U}_s &\coloneqq \left(\begin{matrix} 
    u_s \\ v_s \end{matrix}\right), &
  \vec{U} &\coloneqq \left(\begin{matrix} 
    \vec{U}_f \\ \vec{U}_s \end{matrix}\right),\\
  \boldsymbol{\Phi}_f &\coloneqq \left(\begin{matrix} 
    \varphi_f \\ \psi_f \end{matrix}\right),& 
  \boldsymbol{\Phi}_s &\coloneqq \left(\begin{matrix} 
    \varphi_s \\ \psi_s \end{matrix}\right), &
  \boldsymbol{\Phi} &\coloneqq \left(\begin{matrix} 
    \boldsymbol{\Phi}_f \\ \boldsymbol{\Phi}_s \end{matrix}\right). 
\end{aligned}
\]

After these preliminaries, we are ready to construct a continuous variational formulation of the problem. We define operators describing the fluid and the solid problem
\begin{subequations}
\begin{flalign}
B_f(\vec{U})(\boldsymbol{\Phi}_f) \coloneqq &\int_I \langle \partial_t v_f, \varphi_f \rangle_f \diff t + \int_I a_f(\vec{U})(\boldsymbol{\Phi}_f) \diff t + (v_f(0), \varphi_f(0))_f, \label{b_f} \\
B_s(\vec{U})(\boldsymbol{\Phi}_s) \coloneqq &\int_I \langle \partial_t v_s, \varphi_s \rangle_s \diff t + \int_I \langle \partial_t u_s, \psi_s \rangle_s \diff t + \int_I a_s(\vec{U})(\boldsymbol{\Phi}_s) \diff t \label{b_s} \\
&\qquad + (v_s(0), \varphi_s(0))_s + (u_s(0), \psi_s(0))_s,  \nonumber\\
F_f(\boldsymbol{\Phi}_f) \coloneqq &\int_I ( g_f, \varphi_f )_f \diff t, \nonumber\\
F_s(\boldsymbol{\Phi}_s) \coloneqq &\int_I ( g_s, \varphi_s )_s \diff t \nonumber
\end{flalign}
\end{subequations}
with
\begin{subequations}
\begin{flalign}
a_f(\vec{U})(\boldsymbol{\Phi}_f) & \coloneqq (\nu \nabla v_f, \nabla \varphi_f)_f + (\beta \cdot \nabla v_f, \varphi_f)_f + (\nabla u_f, \nabla \psi_f )_f \label{a_f} \\
&\qquad - \langle \partial_{\vec{n}_f} u_f, \psi_f \rangle_{\Gamma} + \frac{\gamma}{h}\langle u_f - u_s, \psi_f \rangle_{\Gamma} \nonumber \\
&\qquad - \langle \nu \partial_{\vec{n}_f} v_f, \varphi_f \rangle_{\Gamma} + \frac{\gamma}{h} \langle v_f - v_s, \varphi_f \rangle_{\Gamma}, \nonumber \\
a_s(\vec{U})(\boldsymbol{\Phi}_s) & \coloneqq  (\lambda \nabla u_s, \nabla \varphi_s)_s + (\delta \nabla v_s, \nabla \varphi_s)_s - (v_s, \psi_s )_s \label{a_s} \\
& \qquad + \langle \nu \partial_{\vec{n}_f} v_f, \varphi_s \rangle_{\Gamma} - \langle \delta \partial_{\vec{n}_s} v_s, \varphi_s \rangle_{\Gamma}. \nonumber
\end{flalign}
\end{subequations}
All the Laplacian terms were integrated by parts and the dynamic
coupling condition was added. The kinematic coupling condition was
incorporated into the fluid problem, while the dynamic condition
became a part of the solid problem. The Dirichlet boundary conditions
over the interface $\Gamma$ were formulated in a weak sense using
Nitsche's method \cite{Nitsche1971}. We arbitrarily set $\gamma =
1000$, while $h$ is the mesh size.

The compact version of the variational problem presents itself as:
\begin{problem}
Find $\vec{U} \in X$ such that
\begin{flalign*}
& B_f(\vec{U})(\boldsymbol{\Phi}_f) = F_f(\boldsymbol{\Phi}_f) \\
& B_s(\vec{U})(\boldsymbol{\Phi}_s) = F_s(\boldsymbol{\Phi}_s)
\end{flalign*}
for all $\boldsymbol{\Phi}_f \in X_f $ and $\boldsymbol{\Phi}_s \in X_s $. 
\label{continuous_problem}
\end{problem}
%
\subsection{Semi-discrete Petrov-Galerkin formulation}
\label{semi-discrete_formulation}

One of the main challenges emerging from the discretization of Problem \ref{continuous_problem} is the construction of a satisfactory time interval partitioning. Our main objectives include:
\begin{enumerate}
    \item \textbf{Handling coupling conditions} \\
      For the time interval $I = [0, T]$ we introduce a coarse time-mesh which is shared by both of the subproblems
      \[
      0 = t_0 < t_1 < ... < t_N = T,\quad k_n = t_n - t_{n - 1}, \quad I_n = (t_{n - 1}, t_n].
      \]
    We will refer to this mesh as a \textit{macro time mesh}.
    \item \textbf{Allowing for different time-step sizes (possibly non-uniform) in both subproblems} \\
    For each of the subintervals $I_n = (t_{n - 1}, t_{n}]$ we create two distinct submeshes corresponding to each of the subproblems
    $$ t_{n - 1} = t_{f, n}^0 < t_{f, n}^1 < ... < t_{f, n}^{M_n} = t_n, \quad k_{f, n}^m = t_{f, n}^{m} - t_{f, n}^{m - 1}, \quad I_{f, n}^m = (t_{f, n}^{m - 1}, t_{f, n}^{m}],$$
    $$ t_{n - 1} = t_{s, n}^0 < t_{s, n}^1 < ... < t_{s, n}^{L_n} = t_n, \quad k_{s, n}^l = t_{s, n}^{l} - t_{s, n}^{l - 1}, \quad I_{s, n}^l = (t_{s, n}^{l - 1}, t_{s, n}^{l}].$$
    We will refer to these meshes as \textit{micro time meshes}.
\end{enumerate}
We define grid sizes as:
$$k_f : = \max_{n = 1,...,N} \max_{m = 1,...,M_{n}}k_{f, n}^m,\quad
k_s : = \max_{n = 1,...,N} \max_{l = 1,...,L_{n}}k_{s, n}^l,$$
$$k \coloneqq \max \{k_f, k_s \}$$
As trial spaces, we chose spaces consisting of piecewise linear
functions in time,
\[
\begin{aligned}
  X^{1, n}_{f, k}& = \left\{ v \in C(\bar{I_n}, L^2(\Omega_f))|\; v|_{I_{f, n}^m} \in \mathcal{P}_1(I_{f, n}^m, H_f)\text{ for } m = 1,...,M_{n}\right\}, \\
  X^1_{f, k}& = \left\{ v \in C(\bar{I}, L^2(\Omega_f))|\; v|_{I_n} \in X^{1, n}_{f, k} \text{ for } n = 1,...,N\right\}, \\
  X^{1, n}_{s, k}& = \left\{ v \in C(\bar{I_n}, L^2(\Omega_s))|\; v|_{I_{s, n}^l} \in \mathcal{P}_1(I_{s, n}^l, H_s)\text{ for } l = 1,...,L_{n}\right\}, \\
  X^1_{s, k} &= \left\{ v \in C(\bar{I}, L^2(\Omega_s))|\; v|_{I_n} \in X^{1, n}_{s, k} \text{ for } n = 1,...,N \right\},
\end{aligned}
\]
whereas we took spaces of piecewise constant functions as test spaces
\[
\begin{aligned}
  Y^{0, n}_{f, k} &= \left\{ v \in L^2(I_n, L^2(\Omega_f))|\; v|_{I_{f, n}^m} \in \mathcal{P}_0(I_{f, n}^m, H_f) \text{ for } m = 1,...,M_{n} \right. \\
  &\hspace{6cm} \left. \text{ and }v(t_{n - 1}) \in L^2(\Omega_f)\right\}, \\
  Y^{0}_{f, k}& = \left\{ v \in L^2(I, L^2(\Omega_f))|\; v|_{I_{n}} \in Y^{0, n}_{f, k}\text{ for } n = 1,...,N\right\}, \\
  Y^{0, n}_{s, k} &= \left\{ v \in L^2(I_n, L^2(\Omega_s))|\; v|_{I_{s, n}^l} \in \mathcal{P}_0(I_{s, n}^l, H_s)\text{ for }l = 1,...,L_{n} \right.\\
  &\hspace{6cm} \left.  \text{ and }v(t_{n - 1}) \in L^2(\Omega_s)\right.\}, \\
  Y^{0}_{s, k} &= \left\{ v \in L^2(I, L^2(\Omega_s))|\; v|_{I_{n}} \in Y^{0, n}_{s, k}\text{ for } n = 1,...,N\right\}.
\end{aligned}
\]
By $\mathcal{P}_r(I,H)$ we denote the space of polynomials with degree $r$ and values in $H$.
To shorten the notation, we set
\[
\begin{aligned}
  X_{f, k}^n&\coloneqq \left(X_{f, k}^{1, n} \right)^2,\quad&   X_{s, k}^n &\coloneqq \left(X_{s, k}^{1, n} \right)^2,  \quad& X_k^n &\coloneqq X_{f, k}^n \times X_{s, k}^n, \\
  X_{f, k} &\coloneqq \left(X_{f, k}^1 \right)^2,&  X_{s, k}&\coloneqq \left(X_{s, k}^1 \right)^2,& X_k  & \coloneqq X_{f, k} \times X_{s, k}, \\
  Y_{f, k}^n&\coloneqq \left(Y_{f, k}^{0, n} \right)^2,& Y_{s, k}^n& \coloneqq \left(Y_{s, k}^{0, n} \right)^2,& Y_k^n& \coloneqq Y_{f, k}^n \times Y_{s, k}^n, \\
  Y_{f, k}& \coloneqq \left(Y_{f, k}^0 \right)^2,& Y_{s, k}& \coloneqq \left(Y_{s, k}^0 \right)^2,& Y_k& \coloneqq Y_{f, k} \times Y_{s, k}.
\end{aligned}
\]
We assume that inner points of fluid and solid micro time-meshes do not necessarily coincide, i.~e. for every $n = 1, ..., N$, $m = 1, ..., M_{n} - 1$, $l = 1,...,L_{n} - 1$ we may have $t_{f, n}^{m} \neq t_{s, n}^{l}$. Because of this fact, a function defined on the fluid micro time-mesh can not be directly evaluated in the points of the solid micro time mesh, and vice versa. To solve this problem, we introduce nodal interpolation operators
\[
i_n^f:X^n \to X_f^n \times \mathcal{P}_1(I_n, X_s^n),\quad 
i_n^s:X^n \to \mathcal{P}_1(I_n, X_f^n) \times X_s^n,
\]
where $X^n \coloneqq X\Big|_{I_n}$, $X_n^f \coloneqq X^f\Big|_{I_n}$, $X_n^s \coloneqq X^s\Big|_{I_n}$ and
\begin{equation}  \label{interpolation_operator_primal}
  \begin{aligned}
    i_n^f \vec{U}(t) &\coloneqq \left(\begin{matrix} \vec{U}_f(t) \\ \frac{t_n - t}{k_n}\vec{U}_s(t_{n - 1}) + \frac{t - t_{n - 1}}{k_n}\vec{U}_s(t_n) \end{matrix}\right), \\
    i_n^s \vec{U}(t) &\coloneqq \left(\begin{matrix} \frac{t_n - t}{k_n}\vec{U}_f(t_{n - 1}) + \frac{t - t_{n - 1}}{k_n}\vec{U}_f(t_n) \\ 
      \vec{U}_s(t)\end{matrix}\right). 
  \end{aligned}
\end{equation}
Since the operators $B_f$ and $B_s$ are linear, the resulting scheme is equivalent to the Crank-Nicolson scheme up to the numerical quadrature of $F_f$, see also~\cite{ErikssonEstepHansboJohnson1995,Thomee1997}. Taking trial functions piecewise linear in time $\vec{U}_k \in X_k$ and test functions piecewise constant in time $\boldsymbol{\Phi}_{f, k} \in Y_{f, k}$, $\boldsymbol{\Phi}_{s, k} \in Y_{s, k}$, we can construct operators on every of the macro time-steps $I_n = (t_{n - 1}, t_n]$
  \begin{equation}\label{b_f^n}
    \begin{aligned}
      B_f^n(\vec{U}_k)(\boldsymbol{\Phi}_{f, k}) \coloneqq & \sum_{m = 1}^{M_{n}} \bigg\{ (v_{f, k}(t_{f, n}^m) - v_{f, k}(t_{f, n}^{m - 1}), \varphi_{f, k}(t_{f, n}^{m}))_f \\
      & \qquad + \frac{k_{f, n}^m}{2}a_f(i_n^f \vec{U}_{k}(t_{f, n}^m))(\boldsymbol{\Phi}_{f, k}(t_{f, n}^m))  \\
      & \qquad + \frac{k_{f, n}^m}{2}a_f(i_n^f\vec{U}_{k}(t_{f, n}^{m - 1}))(\boldsymbol{\Phi}_{f, k}(t_{f, n}^m))\bigg\},  
    \end{aligned}
  \end{equation}
  \begin{equation}\label{b_s^n}
    \begin{aligned}
      B_s^n(\vec{U}_k)(\boldsymbol{\Phi}_{s, k}) \coloneqq & \sum_{l = 1}^{L_{n}} \bigg\{ (v_{s, k}(t_{s, n}^{l}) - v_{s, k}(t_{s, n}^{l - 1}), \varphi_{s, k}(t_{s, n}^{l})_s \\
      & \qquad + (u_{s, k}(t_{s, n}^{l}) - u_{s, k}(t_{s, n}^{l - 1}), \psi_{s, k}(t_{s, n}^{l}))_s  \\ 
      & \qquad + \frac{k_{s, n}^l}{2}a_s(i_n^s\vec{U}_{k}(t_{s, n}^l))(\boldsymbol{\Phi}_{s, k}(t_{s, n}^l))  \\ 
      & \qquad + \frac{k_{f, n}^l}{2}a_s(i_n^s\vec{U}_{k}(t_{s, n}^{l - 1}))(\boldsymbol{\Phi}_{s, k}(t_{s, n}^l)) \bigg\},
    \end{aligned}
  \end{equation}
  \begin{equation*}
    \begin{aligned}
      F_f^n(\boldsymbol{\Phi}_{f, k}) \coloneqq & \sum_{m = 1}^{M_{n}} \left(\int_{I_{s, n}^m}g_f(t) \diff t, \varphi_{f, k}(t_{f, n}^m) \right)_f,\\
      F_s^n(\boldsymbol{\Phi}_{s, k}) \coloneqq & \sum_{l = 1}^{L_{n}} \left(\int_{I_{s, n}^l}g_s(t) \diff t, \varphi_{s, k}(t_{s, n}^l) \right)_s 
    \end{aligned}
  \end{equation*}
  Then, the forms on the whole time interval $I= [0, T]$ are just sums of the operators over the subintervals and initial conditions:
  \begin{subequations}
    \begin{flalign*}
      B_f(\vec{U}_k)(\boldsymbol{\Phi}_{f, k}) = & \sum_{n = 1}^{N} B_f^n(\vec{U}_k)(\boldsymbol{\Phi}_{f, k}) + (v_{f, k}(t_0), \varphi_{f, k}(t_0))_f,\\
      B_s(\vec{U}_k)(\boldsymbol{\Phi}_{s, k}) = & \sum_{n = 1}^{N}B_s^n(\vec{U}_k)(\boldsymbol{\Phi}_{s, k})
      + (v_{s, k}(t_0), \varphi_{s, k}(t_0))_s + (u_{s, k}(t_0), \psi_{s, k}(t_0))_s, \\
      F_f(\boldsymbol{\Phi}_{f, k}) = & \sum_{n = 1}^{N} F_f^n(\boldsymbol{\Phi}_{f, k}), \\
      F_s(\boldsymbol{\Phi}_{s, k}) = & \sum_{n = 1}^{N} F_s^n(\boldsymbol{\Phi}_{s, k})
    \end{flalign*}
  \end{subequations}
  With that at hand, we can pose a semi-discrete variational problem:
  \begin{problem}
Find $\vec{U}_k \in X_k$ such that:
\begin{flalign*}
& B_f(\vec{U}_k)(\boldsymbol{\Phi}_{f, k}) = F_f(\boldsymbol{\Phi}_{f, k}) \\
& B_s(\vec{U}_k)(\boldsymbol{\Phi}_{s, k}) = F_s(\boldsymbol{\Phi}_{s, k})
\end{flalign*}
for all $\boldsymbol{\Phi}_{f, k} \in Y_{f, k}$ and $\boldsymbol{\Phi}_{s, k} \in Y_{s, k}$. 
\label{semi_discrete_problem}
\end{problem}
%
%
\section{Decoupling methods}
\label{decoupling_methods}
Even though Problem \ref{semi_discrete_problem} is discretized in time, it is still coupled across the interface. That makes solving the subproblems independently impossible. To deal with this obstacle, we chose to use an iterative approach on each of the subintervals $I_n$ and introduce decoupling strategies. For a fixed time interval $I_n$ every iteration of a decoupling method consists of the following steps:
\begin{enumerate}
\item Using the solution of the solid subproblem from the previous iteration $\vec{U}_{s, k}^{(i - 1)}$, we set the boundary conditions on the interface at the time $t_n$, solve the fluid problem and get the solution $\vec{U}_{f, k}^{(i)}$.
\item Similarly, we use the solution $\vec{U}_{f, k}^{(i)}$ for setting the boundary conditions of the solid problem and obtain an intermediate solution $\widetilde{\vec{U}}_{s, k}^{(i)}$.
\item We apply a decoupling function to the intermediate solution $\widetilde{\vec{U}}_{s, k}^{(i)}$ and acquire $\vec{U}_{s, k}^{(i)}$.
\end{enumerate}
This procedure is visualized by
\[ 
\vec{U}_{s, k}^{(i - 1)} \xrightarrow[\text{subproblem}]{\text{fluid}}
\vec{U}_{f, k}^{(i)}
\xrightarrow[\text{subproblem}]{\text{solid}}
\widetilde{\vec{U}}_{s, k}^{(i)}
\xrightarrow[\text{function}]{\text{decoupling}}
\vec{U}_{s, k}^{(i)}.
\]
The main challenge emerges from the transition between $\widetilde{\vec{U}}_{s, k}^{(i)}$ and $\vec{U}_{s, k}^{(i)}$. In the next subsections, we will present two techniques. The first one is the relaxation method described in Section \ref{relaxation}. The second one, in Section \ref{shooting}, is the shooting method.

We clarify how the intermediate solution $\widetilde{\vec{U}}_{s,k}^{(i)}$ is obtained from $\vec{U}_{s, k}^{(i - 1)}$ by the definition of Problem~\ref{decoupled_problem}.
\begin{problem}
  For a given $\vec{U}_{s, k}^{(i - 1)} \in X_{s, k}^n$, find $\vec{U}_{f, k}^{(i)} \in X_{f, k}^n$ and $\widetilde{\vec{U}}_{s, k}^{(i)} \in X_{s, k}^n$ such that:
  \begin{flalign*}
    B_f^n &\left(\begin{array}{l} 
      \vec{U}_{f, k}^{(i)} \\ \vec{U}_{s, k}^{(i - 1)} \end{array}\right)(\boldsymbol{\Phi}_{f, k}) = F_f^n(\boldsymbol{\Phi}_{f, k}) \\
    B_s^n & \left(\begin{array}{l} 
      \vec{U}_{f, k}^{(i)} \\ \widetilde{\vec{U}}_{s, k}^{(i)} \end{array}\right)(\boldsymbol{\Phi}_{s, k}) = F_s^n(\boldsymbol{\Phi}_{s, k})
  \end{flalign*}
  for all $\boldsymbol{\Phi}_{f, k} \in Y_{f, k}^n$ and $\boldsymbol{\Phi}_{s, k} \in Y_{s, k}^n$.
  \label{decoupled_problem}
\end{problem}
\begin{remark}
  \normalfont
  Even though in Problem \ref{decoupled_problem} we demand $\vec{U}_{s, k}^{(i - 1)} \in X_{s, k}^n$, in fact, assuming we already know $\vec{U}_{s, k}(t_{n - 1})$,  it is sufficient to set $\left(\vec{U}_{s,k}^{(i - 1)}(t_n)\right)\Big|_{\Gamma}$. The semi-discrete fluid operator (\ref{b_f^n}) is coupled with the solid operator (\ref{b_s^n}) only across the interface~$\Gamma$. Additionally, the interpolation operator (\ref{interpolation_operator_primal}) constructs values over the whole time interval $I_n$ based only on values in the points $t_{n - 1}$ and $t_n$.
  \label{boundary_values}
\end{remark}


\subsection{Relaxation method}
\label{relaxation}

The first of the presented methods consists of a simple interpolation operator being an example of a fixed point method. It contains the iterated solution of each of the two subproblems, taking the interface values from the last iteration of the other problem. For reasons of stability, such explicit partitioned iteration usually requires the introduction of a damping parameter. Here, we only consider fixed damping parameters. 

\begin{definition}[Relaxation Function]
  Let $\vec{U}_{s, k}^{(i -1)} \in X_{s, k}^n$ and  $\widetilde{\vec{U}}_{s, k}^{(i)} \in X_{s, k}^n$ be the solid solution of Problem \ref{decoupled_problem}. Then for $\tau \in [0, 1]$ the relaxation function $R: X_{s, k}^n \to X_{s, k}^n$ is defined as:
  \begin{equation*}
    R(\vec{U}_{s, k}^{(i -1)})\coloneqq \tau \widetilde{\vec{U}}_{s, k}^{(i)} + (1 - \tau)\vec{U}_{s, k}^{(i -1)}
  \end{equation*}
\end{definition}
Assuming that we already know the value $\vec{U}_{s, k}(t_{n -1})$, we pose
\begin{equation*}\left\{
  \begin{aligned}
    \vec{U}_{s, k}^{(0)}(t_n)&\coloneqq \vec{U}_{s, k}(t_{n - 1}), \\
    \vec{U}_{s, k}^{(i)}(t_n)&\coloneqq R(\vec{U}_{s,k}^{(i-1)})(t_n).
  \end{aligned}\right.
\end{equation*}
The stopping criterion is based on checking how far the computed solution is from the fixed point. We evaluate the $l^{\infty}$ norm of $\left( \widetilde{\vec{U}}_{s, k}^{(i + 1)}(t_n) - \vec{U}_{s, k}^{(i)}(t_n) \right)\Big|_{\Gamma}$ and once for $i_{\text{stop}}$ this norm is desirably small, we set

$$\vec{U}_{k}(t_n) \coloneqq \left(\begin{matrix} 
  \vec{U}_{f, k}^{(i_{\text{stop}})} \\ \vec{U}_{s, k}^{(i_{\text{stop}})} \end{matrix}\right)(t_n).$$

\subsection{Shooting method}
\label{shooting}
Here we present another iterative method, where we define a root-finding problem on the interface. We use the Newton method with a matrix-free GMRES method for approximation of the inverse of the Jacobian.

\begin{definition}[Shooting Function]
  Let $\vec{U}_{s, k}^{(i -1)} \in X_{s, k}^n$ and  $\widetilde{\vec{U}}_{s, k}^{(i)} \in X_{s, k}^n$ be the solid solution of Problem \ref{decoupled_problem}. Then the shooting function $S: (X_{s, k}^n)^2 \to (L^2(\Gamma))^2$ is defined as:
  \begin{equation}
    S(\vec{U}_{s, k}^{(i -1)})\coloneqq \left(\vec{U}_{s,k}^{(i - 1)}(t_n) - \widetilde{\vec{U}}_{s,k}^{(i)}(t_n) \right)\Big|_{\Gamma} \label{shooting_function}
  \end{equation}
\end{definition}

Our aim is finding the root of function (\ref{shooting_function}). To do so, we employ the Netwon method
\begin{equation*}
  S'(\vec{U}_{s,k}^{(i - 1)})\vec{d} =  -S(\vec{U}_{s,k}^{(i - 1)}).
\end{equation*}

In each iteration of the Newton method, the greatest difficulty causes computing and inverting the Jacobian $S'(\vec{U}_{s,k}^{(i - 1)})$. Instead of approximating all entries of the Jacobian matrix, we consider an approximation of the matrix-vector product only. Since the Jacobian matrix-vector product can be interpreted as a directional derivative, one can assume
\begin{equation}
  S'(\vec{U}_{s,k}^{(i - 1)})\vec{d} \approx \frac{S(\vec{U}_{s,k}^{(i - 1)} + \varepsilon \vec{d} ) - S(\vec{U}_{s,k}^{(i - 1)} )}{\varepsilon}.
  \label{jacobian_operator}
\end{equation}
In principle, the vector $\vec{d}$ is not known. Thus, the formula above can not be used for solving the system directly. However, it is possible to use this technique with iterative solvers which only require the computation of matrix-vector products. Because we did not want to assume much structure of the operator (\ref{jacobian_operator}), we chose the matrix-free GMRES method. Such matrix-free Newton-Krylov methods are frequently used if the Jacobian is not available or too costly for evaluation~\cite{KnollKeyes2004}. 
Once $\vec{d}$ is computed, we set 
\begin{equation} \begin{cases}
    \vec{U}_{s, k}^{(0)}(t_n)\big|_{\Gamma}: = \vec{U}_{s, k}(t_{n - 1})\big|_{\Gamma}, \\
    \vec{U}_{s, k}^{(i)}(t_n)\big|_{\Gamma}\coloneqq \vec{U}^{(i - 1)}_{s, k}(t_n)\big|_{\Gamma} + \vec{d}.
  \end{cases}
\end{equation}
Here, we stop iterating when the $l^{\infty}$ norm of $S(\vec{U}_{s, k}^{(i)})$ is sufficiently small and then we take
$$\vec{U}_{k}(t_n)\big|_{\Gamma} \coloneqq \left(\begin{matrix} 
  \vec{U}_{f, k}^{(i_{\text{stop}})} \\ \widetilde{\vec{U}}_{s, k}^{(i_{\text{stop}})} \end{matrix}\right)(t_n)\big|_{\Gamma}.$$
We note that the method presented here is similar to the one presented in \cite{Degroote2009}, where the authors also introduced a root-finding problem on the interface and solved it with a quasi-Newton method. The main difference lies in the approximation of the inverse of the Jacobian. Instead of using a matrix-free linear solver, there the Jacobian is approximated by solving a least-squares problem.


\subsection{Numerical comparison of the performance}
\label{comparison}

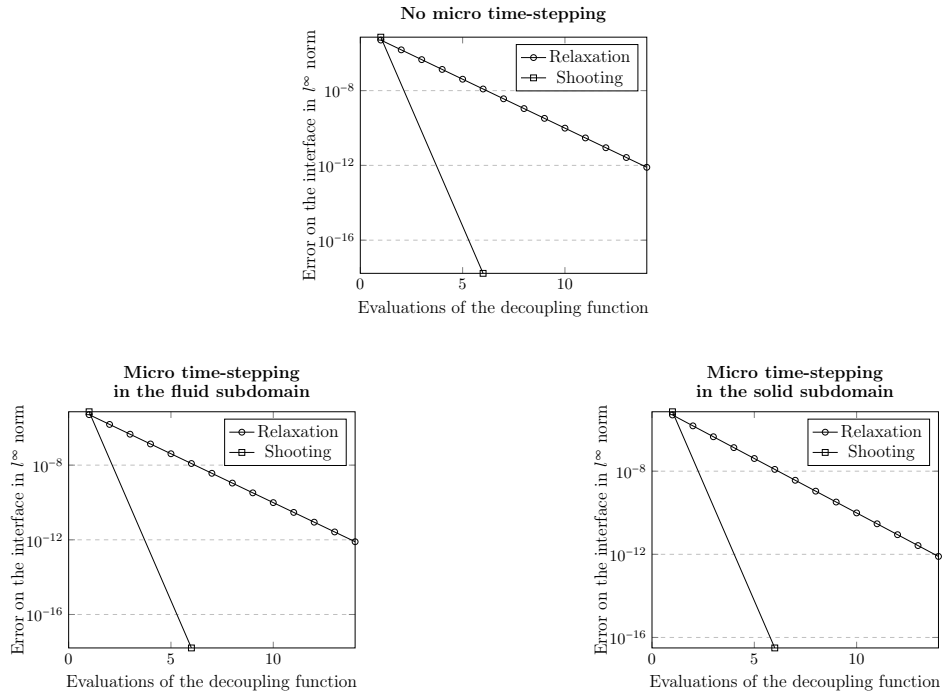
\begin{figure}[t]
  \begin{center}
    \begin{tikzpicture}[scale = 0.55]
      \begin{semilogyaxis}[
	  title style={align=center},
	  title={\large{\textcolor{white}{adjust to the next line}} \\ \large{\textbf{No micro time-stepping}}},
	  xlabel={\large{Evaluations of the decoupling function}},
	  ylabel={\large{Error on the interface in $l^{\infty}$ norm}},
	  xmin=0, xmax=14,
	  ymin=1.6782839494474103e-18, ymax=7.2335262037580485e-06,
	  xtick={0,5,10},
	  ytick={1.0e-16, 1.0e-12, 1.0e-8},
	  legend pos=north east,
	  ymajorgrids=true,
	  grid style=dashed,
	]
	\addplot[
	  color=black,
	  mark=o,
	]
	coordinates {
	  (1.0, 5.063468341658501e-06) 
	  (2.0, 1.5164739991123107e-06) 
	  (3.0, 4.541735599230415e-07) 
	  (4.0, 1.360218687206842e-07) 
	  (5.0, 4.073761844344486e-08) 
	  (6.0, 1.2200638101940217e-08) 
	  (7.0, 3.6540077493352834e-09) 
	  (8.0, 1.09435037943872e-09) 
	  (9.0, 3.2775047633182807e-10) 
	  (10.0, 9.815903496973628e-11) 
	  (11.0, 2.9397962946332215e-11) 
	  (12.0, 8.804490045335847e-12) 
	  (13.0, 2.636884962462999e-12) 
	  (14.0, 7.897291452087574e-13)
	};
	\addplot[
	  color=black,
	  mark=square,
	]
	coordinates {
	  (1.0, 7.2335262037580485e-06) 
	  (6.0, 1.6782839494474103e-18)  
	};
	\legend{\large{Relaxation}, \large{Shooting}}
	
      \end{semilogyaxis}
    \end{tikzpicture}
  \end{center}
  
  \begin{multicols}{2}
    \begin{center}
      
      \begin{tikzpicture}[scale = 0.55]
	\begin{semilogyaxis}[
	    title style={align=center},
	    title={\large{\textbf{Micro time-stepping}}\\\large{\textbf{in the fluid subdomain}}},
	    xlabel={\large{Evaluations of the decoupling function}},
	    ylabel={\large{Error on the interface in $l^{\infty}$ norm}},
	    xmin=0, xmax=14,
	    ymin=1.6143131874927976e-18, ymax=7.2336889009813875e-06,
	    xtick={0,5,10},
	    ytick={1.0e-16, 1.0e-12, 1.0e-8},
	    legend pos=north east,
	    ymajorgrids=true,
	    grid style=dashed,
	  ]
	  \addplot[
	    color=black,
	    mark=o,
	  ]
	  coordinates {
	    (1.0, 5.063582229714357e-06) 
	    (2.0, 1.516508024155365e-06) 
	    (3.0, 4.5418372513637254e-07) 
	    (4.0, 1.3602490562847677e-07) 
	    (5.0, 4.073852572897223e-08) 
	    (6.0, 1.2200909154588915e-08) 
	    (7.0, 3.654088726135697e-09) 
	    (8.0, 1.0943745710675214e-09) 
	    (9.0, 3.277577034094575e-10) 
	    (10.0, 9.816119403036871e-11) 
	    (11.0, 2.939860797104603e-11) 
	    (12.0, 8.804682734489964e-12) 
	    (13.0, 2.6369425288130414e-12) 
	    (14.0, 7.897463161900205e-13)
	  };
	  
	  \addplot[
	    color=black,
	    mark=square,
	  ]
	  coordinates {
	    (1.0, 7.2336889009813875e-06) 
	    (6.0, 1.6143131874927976e-18)  
	  };
	  \legend{\large{Relaxation}, \large{Shooting}}
	  
	\end{semilogyaxis}
      \end{tikzpicture}
    \end{center}
    
    \begin{center}
      \begin{tikzpicture}[scale = 0.55]
	\begin{semilogyaxis}[
	    title style={align=center},
	    title={\large{\textbf{Micro time-stepping}}\\\large{\textbf{in the solid subdomain}}},
	    xlabel={\large{Evaluations of the decoupling function}},
	    ylabel={\large{Error on the interface in $l^{\infty}$ norm}},
	    xmin=0, xmax=14,
	    ymin=3.1299108171965233e-17, ymax=7.234364854672545e-06,
	    xtick={0,5,10},
	    ytick={1.0e-16, 1.0e-12, 1.0e-8},
	    legend pos=north east,
	    ymajorgrids=true,
	    grid style=dashed,
	  ]

	  \addplot[
	    color=black,
	    mark=o,
	  ]
	  coordinates {
	    (1.0, 5.0640553973980905e-06) 
	    (2.0, 1.516686459736926e-06) 
	    (3.0, 4.5424819164532013e-07) 
	    (4.0, 1.3604752339998525e-07) 
	    (5.0, 4.074629352738889e-08) 
	    (6.0, 1.2203533968598318e-08) 
	    (7.0, 3.6549644338047178e-09) 
	    (8.0, 1.0946637381482126e-09) 
	    (9.0, 3.27852382732407e-10) 
	    (10.0, 9.819197494872247e-11) 
	    (11.0, 2.9408554182419644e-11) 
	    (12.0, 8.807880058075925e-12) 
	    (13.0, 2.637965706152789e-12) 
	    (14.0, 7.900728285347777e-13)
	  };
	  
	  \addplot[
	    color=black,
	    mark=square,
	  ]
	  coordinates {
	    (1.0, 7.234364854672545e-06) 
	    (6.0,  3.1299108171965233e-17)  
	  };
	  \legend{\large{Relaxation}, \large{Shooting}}
	  
	\end{semilogyaxis}
      \end{tikzpicture}
    \end{center}
  \end{multicols}
  \caption{Performance of decoupling methods for Configuration \ref{configuration_1} in one macro time-step in the case of $M_{n} = 1$ and $L_{n} = 1$ (top), $M_{n} = 10$ and $L_{ n} = 1$ (left), $M_{n} = 1$ and $L_{n} = 10$ (right).}
  \label{comparison_one_timestep_configuration_1}
\end{figure}

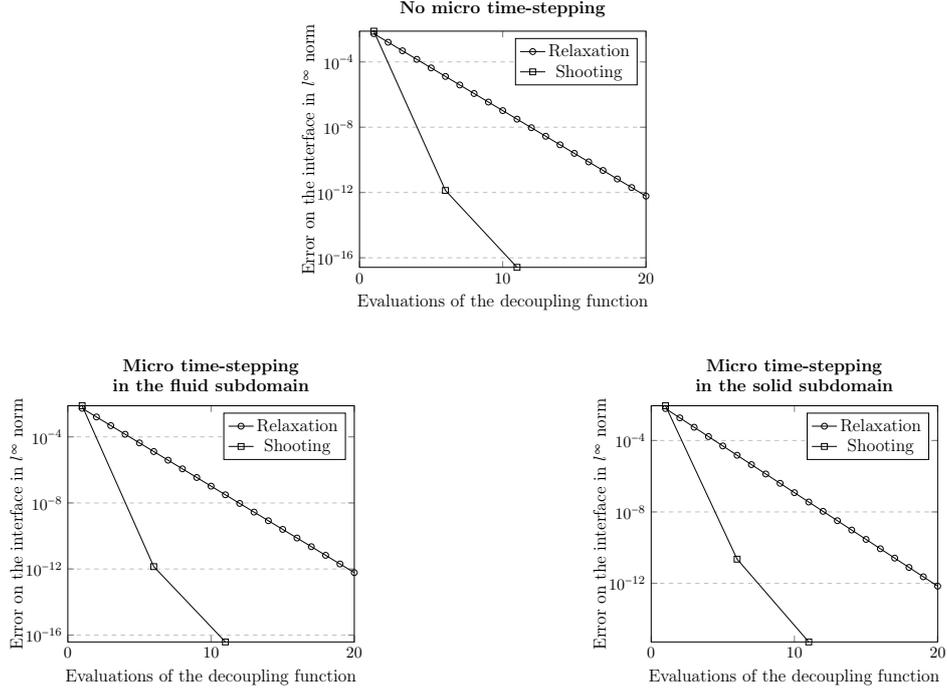
\begin{figure}[t]
  \begin{center}
    \begin{tikzpicture}[scale = 0.55]
      \begin{semilogyaxis}[
	  title style={align=center},
	  title={\large{\textcolor{white}{adjust to the next line}} \\ \large{\textbf{No micro time-stepping}}},
	  xlabel={\large{Evaluations of the decoupling function}},
	  ylabel={\large{Error on the interface in $l^{\infty}$ norm}},
	  xmin=0, xmax=20,
	  ymin=2.5981947607924254e-17, 
	  ymax=0.007656308303176444,
	  xtick={0, 10, 20},
	  ytick={1.0e-16, 1.0e-12, 1.0e-8, 1.0e-4},
	  legend pos=north east,
	  ymajorgrids=true,
	  grid style=dashed,
	]

	\addplot[
	  color=black,
	  mark=o,
	]
	coordinates {
          (1.0,  0.005359415811535416)
          (2.0, 0.0016050925565451867)
          (3.0, 0.00048070951497492957)
          (4.0, 0.00014396779939445696)
          (5.0, 4.311694922066177e-05)
          (6.0, 1.2913105301192662e-05)
          (7.0, 3.8673490456006905e-06)
          (8.0, 1.1582333328392495e-06)
          (9.0, 3.468795986062422e-07)
          (10.0, 1.0388706259737299e-07)
          (11.0, 3.111316455267489e-08)
          (12.0, 9.318090293869872e-09)
          (13.0, 2.790677516931567e-09)
          (14.0, 8.35780819315297e-10)
          (15.0, 2.503082492867421e-10)
          (16.0, 7.496490104927928e-11)
          (17.0, 2.2451248778689066e-11)
          (18.0, 6.723954703037041e-12)
          (19.0, 2.0137587739115056e-12)
          (20.0, 6.030839615136963e-13)
	};
	
	\addplot[
	  color=black,
	  mark=square,
	]
	coordinates {
	  (1.0, 0.007656308303176444)
          (6.0, 1.3362380439487633e-12)
          (11.0, 2.5981947607924254e-17)
	};
	
	\legend{\large{Relaxation}, \large{Shooting}}
	
      \end{semilogyaxis}
    \end{tikzpicture}
  \end{center}

  \begin{multicols}{2}
    \begin{center}
      \begin{tikzpicture}[scale = 0.55]
	\begin{semilogyaxis}[
	    title style={align=center},
	    title={\large{\textbf{Micro time-stepping}}\\\large{\textbf{in the fluid subdomain}}},
	    xlabel={\large{Evaluations of the decoupling function}},
	    ylabel={\large{Error on the interface in $l^{\infty}$ norm}},
	    xmin=0, xmax=20,
	    ymin=3.8236446003307725e-17, 
	    ymax=0.007656309743571194,
	    xtick={0, 10, 20},
	    ytick={1.0e-16, 1.0e-12, 1.0e-8, 1.0e-4},
	    legend pos=north east,
	    ymajorgrids=true,
	    grid style=dashed,
	  ]

	  \addplot[
	    color=black,
	    mark=o,
	  ]
	  coordinates {
            (1.0, 0.005359416820499795)
            (2.0, 0.0016050927552661936)
            (3.0, 0.0004807095435034251)
            (4.0, 0.00014396779865734116)
            (5.0, 4.3116946220042164e-05)
            (6.0, 1.2913103569907667e-05)
            (7.0, 3.867348277715173e-06)
            (8.0, 1.1582330281757546e-06)
            (9.0, 3.468794850001077e-07)
            (10.0, 1.0388702187616845e-07)
            (11.0, 3.111315033098777e-08)
            (12.0, 9.318085433252181e-09)
            (13.0, 2.7906758946635684e-09)
            (14.0, 8.357802788856557e-10)
            (15.0, 2.503080787801899e-10)
            (16.0, 7.4964826958197e-11)
            (17.0, 2.245124566513367e-11)
            (18.0, 6.723933459764385e-12)
            (19.0, 2.0137405319378455e-12)
            (20.0, 6.030960440037066e-13)
	  };
	  
	  \addplot[
	    color=black,
	    mark=square,
	  ]
	  coordinates {
	    (1.0, 0.007656309743571194)
            (6.0, 1.3908079328133397e-12)
            (11.0, 3.8236446003307725e-17) 
	  };
	  \legend{\large{Relaxation}, \large{Shooting}}
	  
	\end{semilogyaxis}
      \end{tikzpicture}
    \end{center}

    \begin{center}
      \begin{tikzpicture}[scale = 0.55]
	\begin{semilogyaxis}[
	    title style={align=center},
	    title={\large{\textbf{Micro time-stepping}}\\\large{\textbf{in the solid subdomain}}},
	    xlabel={\large{Evaluations of the decoupling function}},
	    ylabel={\large{Error on the interface in $l^{\infty}$ norm}},
	    xmin=0, xmax=20,
	    ymin=5.225568472448191e-16, 
	    ymax=0.008832788810731775,
	    xtick={0, 10, 20},
	    ytick={1.0e-16, 1.0e-12, 1.0e-8, 1.0e-4},
	    legend pos=north east,
	    ymajorgrids=true,
	    grid style=dashed,
	  ]

	  \addplot[
	    color=black,
	    mark=o,
	  ]
	  coordinates {
	    (1.0, 0.006182952167512061)
            (2.0, 0.001851894094775634)
            (3.0, 0.0005546722347904078)
            (4.0, 0.00016613331350273555)
            (5.0, 4.9759619095894025e-05)
            (6.0, 1.4903812813475945e-05)
            (7.0, 4.463933820745542e-06)
            (8.0, 1.3370206885870224e-06)
            (9.0, 4.004594325296193e-07)
            (10.0, 1.1994411584002287e-07)
            (11.0, 3.5925215834573944e-08)
            (12.0, 1.0760187730531687e-08)
            (13.0, 3.2228516215462864e-09)
            (14.0, 9.652967270027322e-10)
            (15.0, 2.8912230173214463e-10)
            (16.0, 8.659680184574186e-11)
            (17.0, 2.5937321917881572e-11)
            (18.0, 7.768221735936033e-12)
            (19.0, 2.3266147157063255e-12)
            (20.0, 6.973162229222618e-13)
	  };
	  
	  \addplot[
	    color=black,
	    mark=square,
	  ]
	  coordinates {
	    (1.0, 0.008832788810731775)
            (6.0, 2.2599031490117177e-11)
            (11.0, 5.225568472448191e-16)  
	  };
	  \legend{\large{Relaxation}, \large{Shooting}}
	  
	\end{semilogyaxis}
      \end{tikzpicture}
    \end{center}
  \end{multicols}
  \caption{Performance of decoupling methods for Configuration \ref{configuration_2} in one macro time-step in the case of $M_{n} = 1$ and $L_{n} = 1$ (top), $M_{n} = 10$ and $L_{ n} = 1$ (left), $M_{n} = 1$ and $L_{n} = 10$ (right).}
  \label{comparison_one_timestep_configuration_2}
\end{figure}

In Figures~\ref{comparison_one_timestep_configuration_1} and \ref{comparison_one_timestep_configuration_2} we present the comparison of the performance of both methods based on the number of micro time-steps. We assumed that the micro time-steps have a uniform size. We performed the simulations in the case of no micro time-stepping ($L_n = 1$, $M_n = 1$), micro time-stepping in the fluid subdomain ($M_n = 10$, $L_n = 1$) and the solid subdomain ($M_n = 1$, $L_n = 10$). Figure~\ref{comparison_one_timestep_configuration_1} shows results for the right hand side according to Configuration~\ref{configuration_1}. Figure~\ref{comparison_one_timestep_configuration_2} corresponds to Configuration~\ref{configuration_2}. We investigated one macro time-step $I_2 = [0.02, 0.04]$. We set the relaxation parameter to $\tau = 0.7$. Both methods are very robust concerning the number of micro time-steps. The relaxation method, as expected, has a linear convergence rate. In both cases, despite the nested GMRES method, the performance of the shooting method is much better. For Configuration~\ref{configuration_1}, the relaxation method needs 13 iterations to converge. The shooting method needs only 2 iterations of the Newton method (which is the reason why each of the graphs in Figure~\ref{comparison_one_timestep_configuration_1} displays only two evaluations of the error) and overall requires 6 evaluations of the decoupling function. In the case of Configuration~\ref{configuration_2}, both methods need more iterations to reach the same level of accuracy. The number of iterations of the relaxation method increases to 20 while the shooting method needs 3 iterations of the Newton method and 11 evaluations of the decoupling function. 
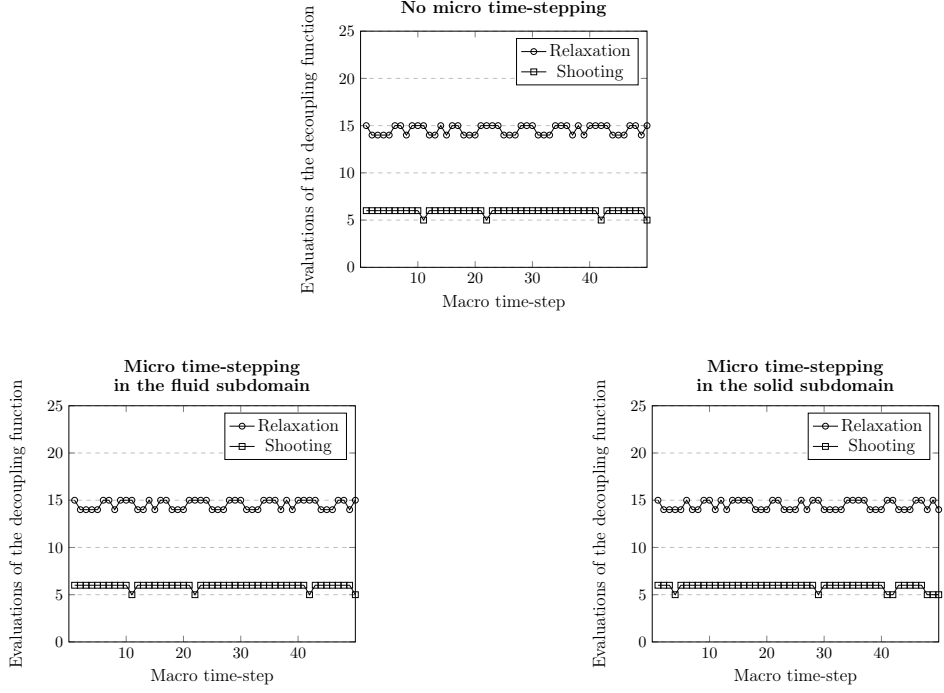
\begin{figure}[t]
  \begin{center}
    \begin{tikzpicture}[scale = 0.55]
      \begin{axis}[
	  title style={align=center},
	  title={\large{\textcolor{white}{adjust to the next line}} \\ \large{\textbf{No micro time-stepping}}},
	  xlabel={\large{Macro time-step}},
	  ylabel={\large{Evaluations of the decoupling function}},
	  xmin=0, xmax=50,
	  ymin=0, ymax=25,
	  xtick={10, 20, 30, 40},
	  ytick={0, 5, 10, 15, 20, 25},
	  legend pos=north east,
	  ymajorgrids=true,
	  grid style=dashed,
	]
        \addplot[
	  color=black,
	  mark=o,
	]
	coordinates {
	  (1, 15) (2, 14) (3, 14) (4, 14) (5, 14) (6, 15) (7, 15) (8, 14) (9, 15) (10, 15) (11, 15) (12, 14) (13, 14) (14, 15) (15, 14) (16, 15) (17, 15) (18, 14) (19, 14) (20, 14) (21, 15) (22, 15) (23, 15) (24, 15) (25, 14) (26, 14) (27, 14) (28, 15) (29, 15) (30, 15) (31, 14) (32, 14) (33, 14) (34, 15) (35, 15) (36, 15) (37, 14) (38, 15) (39, 14) (40, 15) (41, 15) (42, 15) (43, 15) (44, 14) (45, 14) (46, 14) (47, 15) (48, 15) (49, 14) (50, 15)
	};
	
	\addplot[
	  color=black,
	  mark=square,
	]
	coordinates {
	  (1, 6) (2, 6) (3, 6) (4, 6) (5, 6) (6, 6) (7, 6) (8, 6) (9, 6) (10, 6) (11, 5) (12, 6) (13, 6) (14, 6) (15, 6) (16, 6) (17, 6) (18, 6) (19, 6) (20, 6) (21, 6) (22, 5) (23, 6) (24, 6) (25, 6) (26, 6) (27, 6) (28, 6) (29, 6) (30, 6) (31, 6) (32, 6) (33, 6) (34, 6) (35, 6) (36, 6) (37, 6) (38, 6) (39, 6) (40, 6) (41, 6) (42, 5) (43, 6) (44, 6) (45, 6) (46, 6) (47, 6) (48, 6) (49, 6) (50, 5)
	};
	\legend{\large{Relaxation}, \large{Shooting}}
	
      \end{axis}
    \end{tikzpicture}
  \end{center}
  
  \begin{multicols}{2}
    \begin{center}
      \begin{tikzpicture}[scale = 0.55]
	\begin{axis}[
	    title style={align=center},
	    title={\large{\textbf{Micro time-stepping}}\\\large{\textbf{in the fluid subdomain}}},
	    xlabel={\large{Macro time-step}},
	    ylabel={\large{Evaluations of the decoupling function}},
	    xmin=0, xmax=50,
	    ymin=0, ymax=25,
	    xtick={10, 20, 30, 40},
	    ytick={0, 5, 10, 15, 20, 25},
	    legend pos=north east,
	    ymajorgrids=true,
	    grid style=dashed,
	  ]

	  \addplot[
	    color=black,
	    mark=o,
	  ]
	  coordinates {
	    (1, 15) (2, 14) (3, 14) (4, 14) (5, 14) (6, 15) (7, 15) (8, 14) (9, 15) (10, 15) (11, 15) (12, 14) (13, 14) (14, 15) (15, 14) (16, 15) (17, 15) (18, 14) (19, 14) (20, 14) (21, 15) (22, 15) (23, 15) (24, 15) (25, 14) (26, 14) (27, 14) (28, 15) (29, 15) (30, 15) (31, 14) (32, 14) (33, 14) (34, 15) (35, 15) (36, 15) (37, 14) (38, 15) (39, 14) (40, 15) (41, 15) (42, 15) (43, 15) (44, 14) (45, 14) (46, 14) (47, 15) (48, 15) (49, 14) (50, 15)
	  };
	  
	  \addplot[
	    color=black,
	    mark=square,
	  ]
	  coordinates {
	    (1, 6) (2, 6) (3, 6) (4, 6) (5, 6) (6, 6) (7, 6) (8, 6) (9, 6) (10, 6) (11, 5) (12, 6) (13, 6) (14, 6) (15, 6) (16, 6) (17, 6) (18, 6) (19, 6) (20, 6) (21, 6) (22, 5) (23, 6) (24, 6) (25, 6) (26, 6) (27, 6) (28, 6) (29, 6) (30, 6) (31, 6) (32, 6) (33, 6) (34, 6) (35, 6) (36, 6) (37, 6) (38, 6) (39, 6) (40, 6) (41, 6) (42, 5) (43, 6) (44, 6) (45, 6) (46, 6) (47, 6) (48, 6) (49, 6) (50, 5)
	  };
	  \legend{\large{Relaxation}, \large{Shooting}}
	  
	\end{axis}
      \end{tikzpicture}
    \end{center}

    \begin{center}
      \begin{tikzpicture}[scale = 0.55]
	\begin{axis}[
	    title style={align=center},
	    title={\large{\textbf{Micro time-stepping}}\\\large{\textbf{in the solid subdomain}}},
	    xlabel={\large{Macro time-step}},
	    ylabel={\large{Evaluations of the decoupling function}},
	    xmin=0, xmax=50,
	    ymin=0, ymax=25,
	    xtick={10, 20, 30, 40},
	    ytick={0, 5, 10, 15, 20, 25},
	    legend pos=north east,
	    ymajorgrids=true,
	    grid style=dashed,
	  ]

	  \addplot[
	    color=black,
	    mark=o,
	  ]
	  coordinates {
	    (1, 15) (2, 14) (3, 14) (4, 14) (5, 14) (6, 15) (7, 14) (8, 14) (9, 15) (10, 15) (11, 14) (12, 15) (13, 14) (14, 15) (15, 15) (16, 15) (17, 15) (18, 14) (19, 14) (20, 14) (21, 15) (22, 15) (23, 14) (24, 14) (25, 14) (26, 15) (27, 14) (28, 15) (29, 15) (30, 14) (31, 14) (32, 14) (33, 14) (34, 15) (35, 15) (36, 15) (37, 15) (38, 14) (39, 14) (40, 14) (41, 15) (42, 15) (43, 14) (44, 14) (45, 14) (46, 15) (47, 15) (48, 14) (49, 15) (50, 14)
	  };
	  
	  \addplot[
	    color=black,
	    mark=square,
	  ]
	  coordinates {
	    (1, 6) (2, 6) (3, 6) (4, 5) (5, 6) (6, 6) (7, 6) (8, 6) (9, 6) (10, 6) (11, 6) (12, 6) (13, 6) (14, 6) (15, 6) (16, 6) (17, 6) (18, 6) (19, 6) (20, 6) (21, 6) (22, 6) (23, 6) (24, 6) (25, 6) (26, 6) (27, 6) (28, 6) (29, 5) (30, 6) (31, 6) (32, 6) (33, 6) (34, 6) (35, 6) (36, 6) (37, 6) (38, 6) (39, 6) (40, 6) (41, 5) (42, 5) (43, 6) (44, 6) (45, 6) (46, 6) (47, 6) (48, 5) (49, 5) (50, 5)
	  };
	  \legend{\large{Relaxation}, \large{Shooting}}
	  
	\end{axis}
      \end{tikzpicture}
    \end{center}
  \end{multicols}
  \caption{Number of evaluations of the decoupling functions for Configuration \ref{configuration_1} needed for convergence on the time interval $I = [0, 1]$ for $N = 50$ in the case of $M_{n} = 1$ and $L_{n} = 1$ (top), $M_{n} = 10$ and $L_{ n} = 1$ (left), $M_{n} = 1$ and $L_{n} = 10$ (right).}
  \label{comparison_whole_timeline_configuration_1}
\end{figure}
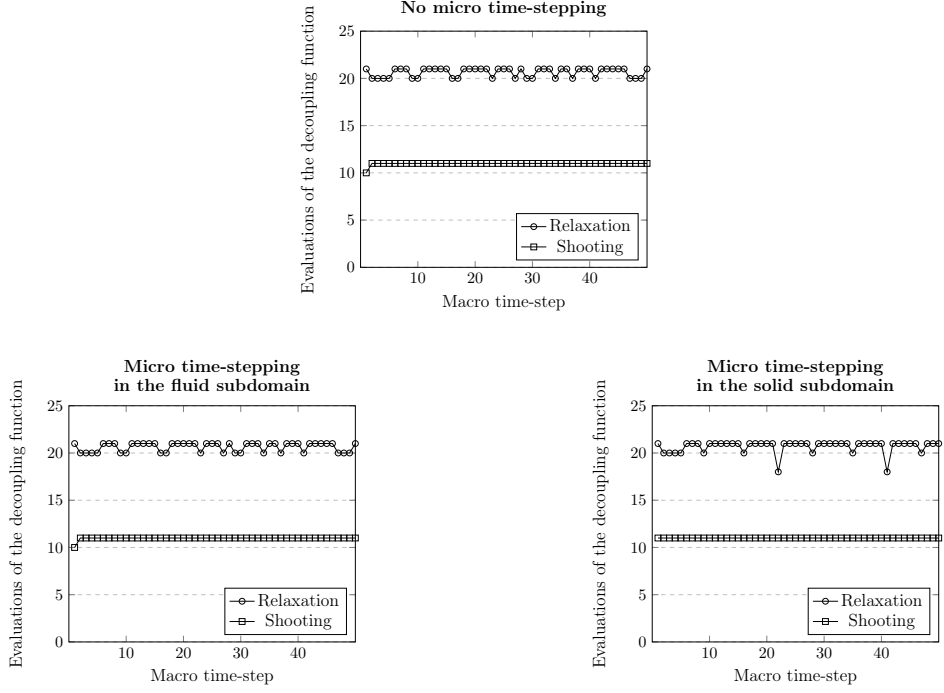
\begin{figure}[t]
  \begin{center}
    \begin{tikzpicture}[scale = 0.55]
      \begin{axis}[
	  title style={align=center},
	  title={\large{\textcolor{white}{adjust to the next line}} \\ \large{\textbf{No micro time-stepping}}},
	  xlabel={\large{Macro time-step}},
	  ylabel={\large{Evaluations of the decoupling function}},
	  xmin=0, xmax=50,
	  ymin=0, ymax=25,
	  xtick={10, 20, 30, 40},
	  ytick={0, 5, 10, 15, 20, 25},
	  legend pos=south east,
	  ymajorgrids=true,
	  grid style=dashed,
	]
        \addplot[
	  color=black,
	  mark=o,
	]
	coordinates {
	  (1, 21) (2, 20) (3, 20) (4, 20) (5, 20) (6, 21) (7, 21) (8, 21) (9, 20) (10, 20) (11, 21) (12, 21) (13, 21) (14, 21) (15, 21) (16, 20) (17, 20) (18, 21) (19, 21) (20, 21) (21, 21) (22, 21) (23, 20) (24, 21) (25, 21) (26, 21) (27, 20) (28, 21) (29, 20) (30, 20) (31, 21) (32, 21) (33, 21) (34, 20) (35, 21) (36, 21) (37, 20) (38, 21) (39, 21) (40, 21) (41, 20) (42, 21) (43, 21) (44, 21) (45, 21) (46, 21) (47, 20) (48, 20) (49, 20) (50, 21)
	};
	\addplot[
	  color=black,
	  mark=square,
	]
	coordinates {
	  (1, 10) (2, 11) (3, 11) (4, 11) (5, 11) (6, 11) (7, 11) (8, 11) (9, 11) (10, 11) (11, 11) (12, 11) (13, 11) (14, 11) (15, 11) (16, 11) (17, 11) (18, 11) (19, 11) (20, 11) (21, 11) (22, 11) (23, 11) (24, 11) (25, 11) (26, 11) (27, 11) (28, 11) (29, 11) (30, 11) (31, 11) (32, 11) (33, 11) (34, 11) (35, 11) (36, 11) (37, 11) (38, 11) (39, 11) (40, 11) (41, 11) (42, 11) (43, 11) (44, 11) (45, 11) (46, 11) (47, 11) (48, 11) (49, 11) (50, 11)
	};
	\legend{\large{Relaxation}, \large{Shooting}}
      \end{axis}
    \end{tikzpicture}
  \end{center}
  \begin{multicols}{2}
    \begin{center}
      \begin{tikzpicture}[scale = 0.55]
	\begin{axis}[
	    title style={align=center},
	    title={\large{\textbf{Micro time-stepping}}\\\large{\textbf{in the fluid subdomain}}},
	    xlabel={\large{Macro time-step}},
	    ylabel={\large{Evaluations of the decoupling function}},
	    xmin=0, xmax=50,
	    ymin=0, ymax=25,
	    xtick={10, 20, 30, 40},
	    ytick={0, 5, 10, 15, 20, 25},
	    legend pos=south east,
	    ymajorgrids=true,
	    grid style=dashed,
	  ]
	  \addplot[
	    color=black,
	    mark=o,
	  ]
	  coordinates {
	    (1, 21) (2, 20) (3, 20) (4, 20) (5, 20) (6, 21) (7, 21) (8, 21) (9, 20) (10, 20) (11, 21) (12, 21) (13, 21) (14, 21) (15, 21) (16, 20) (17, 20) (18, 21) (19, 21) (20, 21) (21, 21) (22, 21) (23, 20) (24, 21) (25, 21) (26, 21) (27, 20) (28, 21) (29, 20) (30, 20) (31, 21) (32, 21) (33, 21) (34, 20) (35, 21) (36, 21) (37, 20) (38, 21) (39, 21) (40, 21) (41, 20) (42, 21) (43, 21) (44, 21) (45, 21) (46, 21) (47, 20) (48, 20) (49, 20) (50, 21)
	  };
	  
	  \addplot[
	    color=black,
	    mark=square,
	  ]
	  coordinates {
	    (1, 10) (2, 11) (3, 11) (4, 11) (5, 11) (6, 11) (7, 11) (8, 11) (9, 11) (10, 11) (11, 11) (12, 11) (13, 11) (14, 11) (15, 11) (16, 11) (17, 11) (18, 11) (19, 11) (20, 11) (21, 11) (22, 11) (23, 11) (24, 11) (25, 11) (26, 11) (27, 11) (28, 11) (29, 11) (30, 11) (31, 11) (32, 11) (33, 11) (34, 11) (35, 11) (36, 11) (37, 11) (38, 11) (39, 11) (40, 11) (41, 11) (42, 11) (43, 11) (44, 11) (45, 11) (46, 11) (47, 11) (48, 11) (49, 11) (50, 11)
	  };
	  \legend{\large{Relaxation}, \large{Shooting}}
	  
	\end{axis}
      \end{tikzpicture}
    \end{center}

    \begin{center}
      \begin{tikzpicture}[scale = 0.55]
	\begin{axis}[
	    title style={align=center},
	    title={\large{\textbf{Micro time-stepping}}\\\large{\textbf{in the solid subdomain}}},
	    xlabel={\large{Macro time-step}},
	    ylabel={\large{Evaluations of the decoupling function}},
	    xmin=0, xmax=50,
	    ymin=0, ymax=25,
	    xtick={10, 20, 30, 40},
	    ytick={0, 5, 10, 15, 20, 25},
	    legend pos=south east,
	    ymajorgrids=true,
	    grid style=dashed,
	  ]

	  \addplot[
	    color=black,
	    mark=o,
	  ]
	  coordinates {
	    (1, 21) (2, 20) (3, 20) (4, 20) (5, 20) (6, 21) (7, 21) (8, 21) (9, 20) (10, 21) (11, 21) (12, 21) (13, 21) (14, 21) (15, 21) (16, 20) (17, 21) (18, 21) (19, 21) (20, 21) (21, 21) (22, 18) (23, 21) (24, 21) (25, 21) (26, 21) (27, 21) (28, 20) (29, 21) (30, 21) (31, 21) (32, 21) (33, 21) (34, 21) (35, 20) (36, 21) (37, 21) (38, 21) (39, 21) (40, 21) (41, 18) (42, 21) (43, 21) (44, 21) (45, 21) (46, 21) (47, 20) (48, 21) (49, 21) (50, 21)
	  };
	  
	  \addplot[
	    color=black,
	    mark=square,
	  ]
	  coordinates {
	    (1, 11) (2, 11) (3, 11) (4, 11) (5, 11) (6, 11) (7, 11) (8, 11) (9, 11) (10, 11) (11, 11) (12, 11) (13, 11) (14, 11) (15, 11) (16, 11) (17, 11) (18, 11) (19, 11) (20, 11) (21, 11) (22, 11) (23, 11) (24, 11) (25, 11) (26, 11) (27, 11) (28, 11) (29, 11) (30, 11) (31, 11) (32, 11) (33, 11) (34, 11) (35, 11) (36, 11) (37, 11) (38, 11) (39, 11) (40, 11) (41, 11) (42, 11) (43, 11) (44, 11) (45, 11) (46, 11) (47, 11) (48, 11) (49, 11) (50, 11)
	  };
	  \legend{\large{Relaxation}, \large{Shooting}}
	  
	\end{axis}
      \end{tikzpicture}
    \end{center}
  \end{multicols}
  \caption{Number of evaluations of the decoupling functions for Configuration \ref{configuration_2} needed for convergence on the time interval $I = [0, 1]$ for $N = 50$ in the case of $M_{n} = 1$ and $L_{n} = 1$ (top), $M_{n} = 10$ and $L_{ n} = 1$ (left), $M_{n} = 1$ and $L_{n} = 10$ (right).}
  \label{comparison_whole_timeline_configuration_2}
\end{figure}

In Figures~\ref{comparison_whole_timeline_configuration_1} and \ref{comparison_whole_timeline_configuration_2} we show the number of evaluations of the decoupling function needed to reach the stopping criteria throughout the complete time interval $I = [0, 1]$ for $N = 50$. Similarly, we performed the simulations in the case of no micro time-stepping, micro time-stepping in the fluid and the solid subdomain. We considered both Configuration~\ref{configuration_1} and~\ref{configuration_2}. In the case of Configuration~\ref{configuration_1}, the number of evaluations of the decoupling function using the relaxation method varied between 14 and 15. For the shooting function, this value was mostly equal to 6 with a few exceptions when only 5 evaluations were needed. For Configuration~\ref{configuration_2}, the relaxation method needed between 18 and 21 iterations while for the shooting method it was almost exactly constant to 11. For each configuration, graphs corresponding to no micro time-stepping and micro time-stepping in the fluid subdomain are the same, while introducing micro time-stepping in the solid subdomain resulted in slight variations. For both decoupling methods, the independence of the performance from the number of micro time-steps extends to the whole time interval $I$.
%
%
\section{Goal oriented estimation}
\label{goal_oriented_estimation}

In Section 1 we formulated the semi-discrete problem enabling usage of different time-step sizes in fluid and solid subdomains, whereas in Section 2 we presented methods designed to efficiently solve such problems. However, so far the choice of the step sizes was purely arbitrary. In this section, we are going to present an easily localized error estimator, which can be used as a criterion for the adaptive choice of the time-step size.

For the construction of the error estimator, we used the dual weighted residual (DWR) method~\cite{BeckerRannacher2001}. Given a differentiable functional $J: X \to \mathbb{R}$, our aim is finding a way to approximate $J(\vec{U}) - J(\vec{U}_k)$, where $\vec{U}$ is the solution to Problem~\ref{continuous_problem} and $\vec{U}_k$ is the solution to Problem~\ref{semi_discrete_problem}.
The goal functional $J: X \to \mathbb{R}$ is split into two parts $J_f: X_f \to \mathbb{R}$ and $J_s:X_s \to \mathbb{R}$ which refer to the fluid and solid subdomains, respectively
$$J(\vec{U}): = J_f(\vec{U}_f) + J_s(\vec{U}_s). $$
The DWR method embeds computing the value of $J$ in the optimal control framework - it is equivalent to solving the following optimization problem
\begin{equation*}
  J(\vec{U}) = \min !, \quad  B(\vec{U})(\boldsymbol{\Phi}) = F(\boldsymbol{\Phi}) \textnormal{ for all }\boldsymbol{\Phi} \in X,
\end{equation*}
where
\begin{flalign*}
  B(\vec{U})(\boldsymbol{\Phi}) & \coloneqq B_f(\vec{U})(\boldsymbol{\Phi}_f) + B_s(\vec{U})(\boldsymbol{\Phi}_s), \\
  F(\boldsymbol{\Phi}) & \coloneqq F_f(\boldsymbol{\Phi}_f) + F_s(\boldsymbol{\Phi}_s).
\end{flalign*}
Solving this problem corresponds to finding stationary points of a Lagrangian $\mathcal{L}: X \times (X \oplus Y_k) \to \mathbb{R}$
\begin{equation*}
  \mathcal{L}(\vec{U}, \vec{Z}): = J(\vec{U}) + F(\vec{Z}) - B(\vec{U})(\vec{Z}). 
\end{equation*}
We can not take $X \times X$ as the domain of $\mathcal{L}$ because we operate in a nonconforming set-up, that is $Y_k \notin X$.
Because the form $B$ describes a linear problem, finding stationary points of $\mathcal{L}$ is equivalent to solving the following problem:
\begin{problem}
  For a given $\vec{U} \in X$ being the solution of Problem
  \ref{continuous_problem}, find $\vec{Z} \in X$ such that: 
  \begin{flalign*}
    B(\boldsymbol{\Xi}, \vec{Z}) = J'_{\vec{U}}(\boldsymbol{\Xi})
  \end{flalign*}
  for all $\boldsymbol{\Xi} \in X$.
  \label{lagrangian_problem}
\end{problem}
The solution $\vec{Z}$ is called an \textit{adjoint solution}. By $J'_{\vec U}(\boldsymbol{\Xi})$ we denote the Gateaux derivative of $J(\cdot)$ at $\vec U$ in direction of the test function $\boldsymbol{\Xi}$.

\subsection{Adjoint problem}
\label{adjoint_problem}

\subsubsection{Continuous variational formulation}
\label{adjoint_continuous_variational_formulation}

As the first step in decoupling the Problem~\ref{lagrangian_problem}, we would like to split the form $B$ into forms corresponding to fluid and solid subproblems. However, we can not fully reuse the forms (\ref{a_f}) and (\ref{a_s}) because of the interface terms - the forms have to be sorted regarding test functions. Thus, after defining abbreviations, 
\[
\begin{aligned}
  \boldsymbol{\Xi}_f &\coloneqq \left(\begin{matrix} 
    \xi_f \\ \eta_f \end{matrix}\right), \quad&
  \boldsymbol{\Xi}_s &\coloneqq \left(\begin{matrix} 
    \xi_s \\ \eta_s \end{matrix}\right), \quad&
  \boldsymbol{\Xi} &\coloneqq \left(\begin{matrix} 
    \boldsymbol{\Xi}_f \\ \boldsymbol{\Xi}_s \end{matrix}\right), \\
  \vec{Z}_f &\coloneqq \left(\begin{matrix} 
    z_f \\ y_f \end{matrix}\right), &
  \vec{Z}_s& \coloneqq \left(\begin{matrix} 
    z_s \\ y_s \end{matrix}\right), &
  \vec{Z} &\coloneqq \left(\begin{matrix} 
    \vec{Z}_f \\ \vec{Z}_s \end{matrix}\right)
\end{aligned}
\]
we choose the splitting
\begin{equation*}
  B(\boldsymbol{\Xi})(\vec{Z})\coloneqq \widetilde{B}_f(\boldsymbol{\Xi}_f)(\vec{Z}) + \widetilde{B}_s(\boldsymbol{\Xi}_s)(\vec{Z}),
\end{equation*}
where
\begin{subequations}
  \begin{flalign*}
    \widetilde{B}_f(\boldsymbol{\Xi}_f)(\vec{Z}) \coloneqq & - \int_I \langle \eta_f, \partial_t z_f \rangle_f \diff t + \int_I \widetilde{a}_f(\boldsymbol{\Xi}_f)(\vec{Z}) \diff t + (\eta_f(T), z_f(T))_f, \\
    \widetilde{B}_s(\boldsymbol{\Xi}_s)(\vec{Z}) \coloneqq & - \int_I \langle \eta_s, \partial_t z_s \rangle_s \diff t - \int_I \langle \xi_s, \partial_t y_s \rangle_s \diff t  + \int_I \widetilde{a}_s(\boldsymbol{\Xi}_s)(\vec{Z}) \diff t \\
    & \qquad + (\eta_s(T), z_s(T))_s + (\xi_s(T), y_s(T))_s \nonumber
  \end{flalign*}
\end{subequations}
and 
\begin{subequations}
  \begin{flalign*}
    \widetilde{a}_f(\boldsymbol{\Xi}_f)(\vec{Z})  \coloneqq & \; (\nu \nabla \eta_f, \nabla z_f)_f + (\beta \cdot \nabla \eta_f, z_f)_f + (\nabla \xi_f, \nabla y_f)_f \\
    &\qquad - \langle \partial_{\vec{n}_f} \xi_f, y_f \rangle_{\Gamma} + \frac{\gamma}{h} \langle \xi_f, y_f \rangle_{\Gamma} - \langle \nu \partial_{\vec{n}_f} \eta_f, z_f \rangle_{\Gamma} + \frac{\gamma}{h}\langle \eta_f, z_f \rangle_{\Gamma} \nonumber \\ 
    & \qquad + \langle \nu \partial_{\vec{n}_f} \eta_f, z_s \rangle_{\Gamma}, \nonumber \\
    \widetilde{a}_s(\boldsymbol{\Xi}_s)(\vec{Z})  \coloneqq & \; (\lambda \nabla \xi_s, \nabla z_s)_s + (\delta \nabla \eta_s, \nabla z_s)_s - (\eta_s, y_s)_s \\
    &\qquad  - \frac{\gamma}{h} \langle \xi_s, y_f \rangle_{\Gamma} - \frac{\gamma}{h} \langle \eta_s, z_f \rangle_{\Gamma} - \langle \delta \partial_{\vec{n}_s} \eta_s, z_s \rangle_{\Gamma}. \nonumber
  \end{flalign*}
\end{subequations}
We have applied integration by parts in time which reveals that the
adjoint problem runs backward in time. That leads to the formulation
of a continuous adjoint variational problem:
\begin{problem}
  For a given $\vec{U} \in X$ being the solution of Problem
  \ref{continuous_problem}, find $\vec{Z} \in X$ such that: 
  \[
  \begin{aligned}
    \widetilde{B}_f(\boldsymbol{\Xi}_f)(\vec{Z}) &= (J_f)'_{\vec{U}}(\boldsymbol{\Xi}_f) \\
    \widetilde{B}_s(\boldsymbol{\Xi}_s)(\vec{Z}) &= (J_s)'_{\vec{U}}(\boldsymbol{\Xi}_s)
  \end{aligned}
  \]
  for all $\boldsymbol{\Xi}_f \in X_f $ and $\boldsymbol{\Xi}_s \in X_s $. 
  \label{adjoint_continuous_problem}
\end{problem}

\subsubsection{Semi-discrete Petrov-Galerkin formulation}
\label{adjoint_semi_discrete_petrov_galerkin_formulation}
The semi-discrete formulation for the adjoint problem is similar to the one of the primal problem. The main difference lies in the fact that this time trial functions are piecewise constant in time $\vec{Z}_k \in Y_k$, while test functions are piecewise linear in time $\boldsymbol{\Xi}_f \in X_{f, k}$, $\boldsymbol{\Xi}_s \in X_{s, k}$. After the rearrangement of the terms in accordance to test functions on every interval $I_n$, we arrive with the scheme 
\begin{subequations}
  \begin{flalign*}
    \widetilde{B}_f^n(\boldsymbol{\Xi}_{f, k})(\vec{Z}_k) = & \
    \frac{k_{f, n}^{M_{n}}}{2}\widetilde{a}_f(\boldsymbol{\Xi}_{f, k}(t_n))(i_n^f\vec{Z}_{k}(t_n))  \\
    & \quad+ \sum_{m = 1}^{M_{n} - 1} \bigg\{ (\eta_{f, k}(t^m_{f, n}), z_{f, k}(t^m_{f, n}) - z_{f, k}(t^{m + 1}_{f, n}))_f \nonumber \\
    & \qquad \qquad +\frac{k^m_{f, n}}{2}\widetilde{a}_f(\boldsymbol{\Xi}_{f, k}(t_{f, n}^m))(i_n^f\vec{Z}_{k}(t_{f, n}^{m})) \nonumber\\
    & \qquad \qquad + \frac{k^{m + 1}_{f, n}}{2}\widetilde{a}_f(\boldsymbol{\Xi}_{f, k}(t_{f, n}^m))(i_n^f\vec{Z}_{k}(t_{f, n}^{m + 1})) \bigg\} \nonumber \\
    & \quad+ (\eta_{f, k}(t_{n - 1}), z_{f, k}(t_{n - 1}) - z_{f, k}(t_{f, n}^{1}))_f \nonumber\\
    &\quad+ \frac{k^1_{f, n}}{2}\widetilde{a}_f(\boldsymbol{\Xi}_{f, k}(t_{n - 1}))(i_n^f\vec{Z}_{k}(t_{f, n}^{1})), \nonumber 
    \\ 
    \widetilde{B}^n_s(\boldsymbol{\Xi}_{s, k})(\vec{Z}_k)
    = & \ \frac{k_{s, n}^{L_n}}{2}\widetilde{a}_s(\boldsymbol{\Xi}_{s, k}(t_n))(i_n^s\vec{Z}_{k}(t_n)) \\
    & \quad+ \sum_{l = 1}^{L_n - 1} \bigg\{ (\eta_{s, k}(t_{s, n}^l), z_{s, k}(t_{s, n}^l) - z_{s, k}(t_{s, n}^{l + 1}))_s \nonumber \\
    & \qquad \qquad + (\xi_{s,k}(t_{s, n}^l), y_{s, k}(t_{s, n}^l) - y_{s, k}(t_{s, n}^{l + 1}))_s \nonumber \\
    & \qquad \qquad + \frac{k^l_{s, n}}{2}\widetilde{a}_s(\boldsymbol{\Xi}_{s, k}(t_{s, n}^l))(i_n^s\vec{Z}_{k}(t_{s, n}^l)) \nonumber\\
    & \qquad \qquad + \frac{k^{l + 1}_{s, n}}{2}\widetilde{a}_s(\boldsymbol{\Xi}_{s, k}(t^l_{s, n}))(i_n^s\vec{Z}_{k}(t^{l + 1}_{s, n})) \bigg\} \nonumber \\
    & \quad+ (\eta_{s, k}(t_{n - 1}), z_{s, k}(t_{n - 1}) - z_{s, k}(t^1_{s, n}))_s \nonumber \\
    & \quad + (\xi_{s, k}(t_{n - 1}), y_{s, k}(t_{n - 1}) - y_{s, k}(t_{s, n}^1))_s \nonumber \\
    & \quad + \frac{k^1_{s, n}}{2}\widetilde{a}_s(\boldsymbol{\Xi}_{s, k}(t_{n - 1}))(i_n^s\vec{Z}_{s, k}(t_{s, n}^1)). \nonumber
  \end{flalign*}
\end{subequations}
Note that the adjoint problem does not have a designated initial value at the final time $T$. Instead, the starting value is implicitly defined by the variational formulation. The final schemes are constructed as sums over the macro time intervals $I_n$ and values at the final time $T$
\begin{subequations}
  \begin{flalign*}
    \widetilde{B}_f(\boldsymbol{\Xi}_{f, k})(\vec{Z}_k) = & \sum_{n = 1}^{N} \widetilde{B}_f^n(\boldsymbol{\Xi}_{f, k})(\vec{Z}_k) + (\eta_{f, k}(T), z_{f, k}(T))_f,\\
    \widetilde{B}_s(\boldsymbol{\Xi}_{s, k})(\vec{Z}_{s, k}) = & \sum_{n = 1}^{N}\widetilde{B}_s^n(\boldsymbol{\Xi}_{s, k})(\vec{Z}_{k})
    + (\eta_{s, k}(T), z_{s, k}(T))_s + (\xi_{s, k}(T), y_{s,k}(T))_s. 
  \end{flalign*}
\end{subequations}
With that at our disposal, we can formulate a semi-discrete adjoint variational problem:
\begin{problem}
  For a given $\vec{U} \in X$ being the solution of Problem \ref{continuous_problem}, find $\vec{Z}_k \in Y_k$ such that:
  \[
  \begin{aligned}
    \widetilde{B}_f(\boldsymbol{\Xi}_{f, k})(\vec{Z}_{k}) &= (J_f)'_{\vec{U}}(\boldsymbol{\Xi}_{f, k}) \\
    \widetilde{B}_s(\boldsymbol{\Xi}_{s,k})(Z_k) &= (J_s)'_{\vec{U}}(\boldsymbol{\Xi}_{s,k})
  \end{aligned}
  \]
  for all $\boldsymbol{\Xi}_{f, k} \in X_{f, k}$ and $\boldsymbol{\Xi}_{s, k} \in X_{s, k}$. 
  \label{adjoint_semi_discrete_problem}
\end{problem}
After formulating the problem in a semi-discrete manner, the decoupling methods from Section~\ref{decoupling_methods} can be applied.

\subsection{A posteriori error estimate}
\label{aposteriori_error}
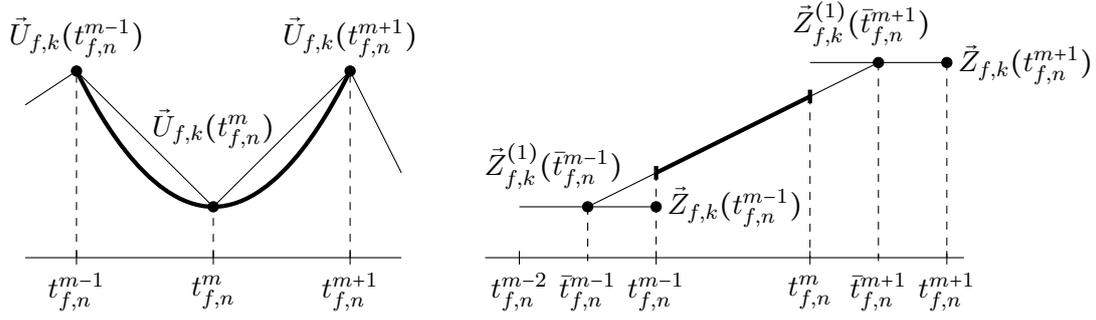
\begin{figure}[t]
  \begin{center}
    \begin{tikzpicture}[scale=0.9]
      \draw (-2.75, -0.75) -- (2.75, -0.75);
      \draw[ultra thick] (0,0) parabola (2,2);
      \draw[ultra thick] (0,0) parabola (-2,2);
      \draw[fill = black] (0,0) circle (0.075cm);
      \draw[fill = black] (2,2) circle (0.075cm);
      \draw[fill = black] (-2,2) circle (0.075cm);
      \draw (-2, -0.65) -- (-2, -0.85);
      \draw (0, -0.65) -- (0, -0.85);
      \draw (2, -0.65) -- (2, -0.85);
      \node at (0, - 1.25) {$t^{m}_{f, n}$};
      \node at (-2, - 1.25) {$t^{m - 1}_{f, n}$};
      \node at (2, - 1.25) {$t^{m + 1}_{f, n}$};
      \draw[dashed] (-2,-0.75) -- (-2,2);
      \draw[dashed] (0,-0.75) -- (0,0);
      \draw[dashed] (2,-0.75) -- (2,2);
      \draw (-2,2) -- (0,0) -- (2,2);
      \draw (-2.75, 1.5) -- (-2, 2);
      \draw (2, 2) -- (2.75, 0.5);
      \node at (2,2.5) {$\vec{U}_{f, k}(t^{m + 1}_{f, n})$};
      \node at (-2,2.5) {$\vec{U}_{f, k}(t^{m - 1}_{f, n})$};
      \node at (0,1.15) {$\vec{U}_{f, k}(t^{m}_{f, n})$};
    \end{tikzpicture}\hspace{0.5cm}
    \begin{tikzpicture}[scale=0.9]
      \draw (0.5, 0.0) -- (7.5, 0.0);
      \draw (1.0, -0.1) -- (1.0, 0.1);
      \draw (3.0, -0.1) -- (3.0, 0.1);
      \draw (5.25, -0.1) -- (5.25, 0.1);
      \draw (7.25, -0.1) -- (7.25, 0.1);
      \draw (2.0, -0.05) -- (2.0, 0.05);
      \draw (6.25, -0.05) -- (6.25, 0.05);
      \node at (1.0, - 0.5) {$  t^{m - 2}_{f, n}$};
      \node at (2.0, - 0.5) {$\bar{t}^{m - 1}_{f, n}$};
      \node at (3.0, - 0.5) {$  t^{m - 1}_{f, n}$};
      \node at (5.25, - 0.5) {$  t^{m}_{f, n}$};
      \node at (6.25, - 0.5) {$\bar{t}^{m + 1}_{f, n}$};
      \node at (7.25, - 0.5) {$t^{m + 1}_{f, n}$};
      \draw[fill = black] (3,0.75) circle (0.075cm);
      \draw[fill = black] (2,0.75) circle (0.075cm);
      \draw[fill = black] (6.25,2.875) circle (0.075cm);
      \draw[fill = black] (7.25,2.875) circle (0.075cm);
      \draw (2, 0.75) -- (3, 1.25);
      \draw (5.25, 2.375) -- (6.25, 2.875);
      \draw[ultra thick] (3, 1.25) -- (5.25, 2.375);
      \draw[ultra thick] (3, 1.15) -- (3, 1.35); 
      \draw[ultra thick] (5.25, 2.275) -- (5.25, 2.475);
      \draw[dashed] (3, 1.25) -- (3, 0);
      \draw (1, 0.75) -- (3, 0.75);
      \draw[dashed] (2, 0.75) -- (2, 0);
      \draw[dashed] (6.25,2.875) -- (6.25, 0);
      \draw[dashed] (7.25,2.875) -- (7.25, 0);
      \draw (5.25,2.875) -- (7.25,2.875);
      \draw[dashed] (5.25, 2.375) -- (5.25, 0);
      \node at (8.4,2.85) {$\vec{Z}_{f, k}(t^{m + 1}_{f, n})$};
      \node at (1.5,1.35) {$\vec{Z}^{(1)}_{f, k}(\bar{t}^{m - 1}_{f, n})$};
      \node at (4.15,0.75) {$\vec{Z}_{f, k}(t^{m - 1}_{f, n})$};
      \node at (6,3.45) {$\vec{Z}^{(1)}_{f, k}(\bar{t}^{m + 1}_{f,n})$};
    \end{tikzpicture}
  \end{center}
  \caption{Reconstruction of the primal solution $\vec{U}_{f, k}^{(2)}$ (left) and the adjoint solution $\vec{Z}_{f, k}^{(1)}$ (right).}
  \label{reconstruction}
\end{figure}

We define the primal residual, split into parts corresponding to the fluid and solid subproblems
\begin{equation*}
  \rho(\vec{U})(\boldsymbol{\Phi}) \coloneqq \rho_f(\vec{U})(\boldsymbol{\Phi}_f) + \rho_s(\vec{U})(\boldsymbol{\Phi}_s), 
\end{equation*}
where
\begin{flalign*}
  \rho_f(\vec{U})(\boldsymbol{\Phi}_f) & \coloneqq F_f(\boldsymbol{\Phi}_f)  - B_f(\vec{U})(\boldsymbol{\Phi}_f), \\
  \rho_s(\vec{U})(\boldsymbol{\Phi}_s) & \coloneqq F_s(\boldsymbol{\Phi}_s) - B_s(\vec{U})(\boldsymbol{\Phi}_s). \nonumber
\end{flalign*}
Similarly, we establish the adjoint residual resulting from the adjoint problem
\begin{equation*}
  \rho^*(\vec{Z})(\boldsymbol{\Xi})\coloneqq \rho_f^*(\vec{Z})(\boldsymbol{\Xi}_f) + \rho_s^*(\vec{Z})(\boldsymbol{\Xi}_s)
\end{equation*}
with
\begin{flalign*}
  \rho_f^*(\vec{Z})(\boldsymbol{\Xi}_f) & \coloneqq (J_f)'_{\vec{U}}(\boldsymbol{\Xi}_f) - \widetilde{B}_f(\boldsymbol{\Xi}_f)(\vec{Z}) \\
  \rho_s^*(\vec{Z})(\boldsymbol{\Xi}_s) & \coloneqq (J_s)'_{\vec{U}}(\boldsymbol{\Xi}_s)- \widetilde{B}_s(\boldsymbol{\Xi}_s)(\vec{Z}). \nonumber
\end{flalign*}
Becker and Rannacher \cite{BeckerRannacher2001} introduced the a posteriori error representation:
\begin{multline}
  J(\vec{U}) - J(\vec{U}_k) = \frac{1}{2} \min_{\boldsymbol{\Phi}_k
    \in Y_{ k}}\rho(\vec{U}_k)(\vec{Z} - \boldsymbol{\Phi}_k) +
  \frac{1}{2}\min_{\boldsymbol{\Xi}_k \in
    X_k}\rho^*(\vec{Z}_k)(\vec{U} - \boldsymbol{\Xi}_k)\\
  + \mathcal{O}(|\vec{U} - \vec{U}_k|^3, |\vec{Z} - \vec{Z}_k|^3) \label{estimator}
\end{multline}
This identity can be used to derive an a posteriori error estimate. Two steps of approximation are required: first, the third order remainder is neglected and second, the approximation errors $\vec Z-\boldsymbol{\Phi}_k$ and $\vec U-\boldsymbol{\Xi}_k$, the \emph{weights}, are replaced by interpolation errors $\vec Z-i_k\vec Z$ and $\vec U-i_k \vec U$, which are then replaced by discrete reconstructions, since the exact solutions $\vec U, \vec Z\in X$ are not available. See \cite{MeidnerRichter2014} and \cite{SchmichVexler2008} for a discussion of different reconstruction schemes. Due to these approximation steps,  this estimator is not precise and it does not result in rigorous bounds. The estimator consists of a primal and adjoint component. Each of them is split again into a fluid and a solid counterpart 
\begin{equation}
  \sigma_k \coloneqq \theta_{f, k} + \theta_{s, k} + \vartheta_{f, k} + \vartheta_{s, k}. \label{residualds_formula}
\end{equation}
The primal estimators are derived from the primal residuals using $\vec{U}_k$ and $\vec{Z}_k$ being the solutions to Problems~\ref{semi_discrete_problem} and \ref{adjoint_semi_discrete_problem}, respectively
\begin{flalign*}
  \theta_{f,k} & \coloneqq \frac{1}{2} \rho_f(\vec{U}_k)(\vec{Z}_{f, k}^{(1)} - \vec{Z}_{f, k}), \\
  \theta_{s,k} & \coloneqq \frac{1}{2} \rho_s(\vec{U}_k)(\vec{Z}_{s, k}^{(1)} - \vec{Z}_{s, k}). \nonumber
\end{flalign*}
The adjoint reconstructions $\vec{Z}_{f, k}^{(1)}$ and $\vec{Z}_{s, k}^{(1)}$ approximating the exact solution are constructed from $\vec{Z}_k$ using linear extrapolation (see Figure~\ref{reconstruction}, right)
\begin{flalign*}
  \vec{Z}_{f, k}^{(1)}\big|_{I_{f, n}^m} \coloneqq & \frac{t - \bar{t}^{m + 1}_{f, n}}{\bar{t}^{m - 1}_{f, n} - \bar{t}^{m + 1}_{f, n}}\vec{Z}_{f, k}(t^{m - 1}_{f, n}) + \frac{t - \bar{t}^{m - 1}_{f, n}}{\bar{t}^{m + 1}_{f, n} - \bar{t}^{m - 1}_{f, n}}\vec{Z}_{f, k}(t^{m + 1}_{f, n}), \\
  \vec{Z}_{s, k}^{(1)}\big|_{I_{s, n}^m} \coloneqq & \frac{t - \bar{t}^{m + 1}_{s, n}}{\bar{t}^{m - 1}_{s, n} - \bar{t}^{m + 1}_{s, n}}\vec{Z}_{s, k}(t^{m - 1}_{s, n}) + \frac{t - \bar{t}^{m - 1}_{s, n}}{\bar{t}^{m + 1}_{s, n} - \bar{t}^{m - 1}_{s, n}}\vec{Z}_{s, k}(t^{m + 1}_{s, n}), \nonumber
\end{flalign*}
with the interval midpoints
\begin{equation}\label{midpoints}
  \bar{t}^m_{f, n} = \frac{t^m_{f, n} + t^{m - 1}_{f, n}}{2},\qquad
  \bar{t}^m_{s, n} = \frac{t^m_{s, n} + t^{m - 1}_{s, n}}{2}.
\end{equation}
The adjoint estimators are based on the adjoint residuals
\begin{flalign*}
  \vartheta_{f,k} \coloneqq \frac{1}{2}\rho_f^*(\vec{Z}_k)(\vec{U}_{f, k}^{(2)} - \vec{U}_{f, k}), \\
  \vartheta_{s,k} \coloneqq \frac{1}{2}\rho_s^*(\vec{Z}_k)(\vec{U}_{s, k}^{(2)} - \vec{U}_{s, k}). \nonumber
\end{flalign*}
The primal reconstructions $\vec{U}_{f, k}^{(2)}$ and $\vec{U}_{s, k}^{(2)}$ are extracted from $\vec{U}_k$ using quadratic reconstruction. The reconstruction is performed on the micro time mesh level on local patches consisting of two neighboring micro time-steps (see Figure~\ref{reconstruction}, left). In general, the patch structure does not have to coincide with the micro and macro time mesh structure - two micro time-steps being in the same local patch do not have to be in the same macro time-step. Additionally, we demand two micro time steps from the same local patch to have the same length.

We compute the effectivity of the error estimate using
\begin{equation*}
  \textnormal{eff}_k \coloneqq \frac{\sigma_k}{J(\vec{U}_{\textnormal{exact}}) - J(\vec{U}_k)} ,
\end{equation*}
where $J(\vec{U}_{\textnormal{exact}})$ can be approximated by extrapolation in time. 

\subsection{Adaptivity}
\label{adaptivity}

The residuals (\ref{residualds_formula}) can be easily localised by restricting them to a specific subinterval
\begin{alignat*}{2}
  & \theta_{f, k}^{n, m}\coloneqq \theta_{f, k}|_{I_{f, n}^{m}}, \qquad && \theta_{s, k}^{n, m}\coloneqq \theta_{s, k}|_{I_{s, n}^{m}}, \\
  & \vartheta_{f, k}^{n, m}\coloneqq \vartheta_{f, k}|_{I_{f, n}^{m}}, \qquad && \vartheta_{s, k}^{n, m}\coloneqq \vartheta_{s, k}|_{I_{s, n}^{m}}.
\end{alignat*}
After defining global numbers of subintervals $M\coloneqq \sum_{n = 1}^N M_n$ and $L\coloneqq \sum_{n = 1}^N L_n$ we can compute an average for each of the components
\begin{equation}
  \bar{\sigma}_{k}\coloneqq \frac{1}{2M} \sum_{n = 1}^{N} \sum_{m = 1}^{M_n} \left( |\theta_{f, k}^{n, m}| + |\vartheta_{f, k}^{n, m}| \right) + \frac{1}{2L} \sum_{n = 1}^{N} \sum_{l = 1}^{L_n} \left( | \theta_{s, k}^{n, l}| + |\vartheta_{s, k}^{n, l}|\right). \label{partial_average}
\end{equation}
This way we can obtain satisfactory refining criteria
\begin{flalign}
  & \left( \left|\theta_{f, k}^{n, m} \right|  \geq  \bar{\sigma}_k \textnormal{ or } \left|\vartheta_{f, k}^{n, m} \right|  \geq  \bar{\sigma}_k \right) \Longrightarrow \textnormal{ refine } I_{f, n}^m, \label{refine_criterium} \\
  & \left( \left|\theta_{s, k}^{n, l} \right|  \geq  \bar{\sigma}_k \textnormal{ or } \left|\vartheta_{s, k}^{n, l} \right|  \geq  \bar{\sigma}_k \right) \Longrightarrow \textnormal{ refine } I_{s, n}^l. \nonumber
\end{flalign}
Taking into account the time interval partitioning structure, we arrive with the following algorithm:
\begin{enumerate}
\item Mark subintervals using the refining criteria (\ref{refine_criterium}).
\item Adjust the local patch structure  - in case only one subinterval from a specific patch is marked, mark the other one as well (see Figure~\ref{local_patches}).
  \begin{figure}[t]
    \centering
    \begin{tikzpicture}[scale = 0.8]
      \draw (-2.75, -0.75) -- (2.75, -0.75);
      \draw[ultra thick] (0,0) parabola (2,2);
      \draw[ultra thick] (0,0) parabola (-2,2);
      \draw[fill = black] (0,0) circle (0.075cm);
      \draw[fill = black] (2,2) circle (0.075cm);
      \draw[fill = black] (-2,2) circle (0.075cm);
      \draw (-2, -0.65) -- (-2, -0.85);
      \draw (0, -0.65) -- (0, -0.85);
      \draw (2, -0.65) -- (2, -0.85);
      \node at (0, - 1.25) {$ t^{m}_{f, n}$};
      \node at (-2, - 1.25) {$t^{m - 1}_{f, n}$};
      \node at (2, - 1.25) {$ t^{m + 2}_{f, n}$};
      \node at (1, - 1.25) {$ t^{m + 1}_{f, n}$};
      \draw[dashed] (-2,-0.75) -- (-2,2);
      \draw[dashed] (0,-0.75) -- (0,0);
      \draw[dashed] (2,-0.75) -- (2,2);
      \draw (1, -0.65) -- (1, -0.85);
      \draw[->] (1.0, -2.1) -- (1.0, -1.65);
      \node at (1.0, -2.35) {refine};
      \draw[thick, ->] (2.75,0.75) .. controls (3.25,1.25) and (4.25,1.25) .. (4.75,0.75);
      
      \draw (4.75, -0.75) -- (9.75, -0.75);
      \draw[ultra thick] (7.5,0) parabola (9.5,2);
      \draw[ultra thick] (7.5,0) parabola (5.5,2);
      \draw[fill = black] (7.5,0) circle (0.075cm);
      \draw[fill = black] (9.5,2) circle (0.075cm);
      \draw[fill = black] (5.5,2) circle (0.075cm);
      \draw (5.5, -0.65) -- (5.5, -0.85);
      \draw (7.5, -0.65) -- (7.5, -0.85);
      \draw (9.5, -0.65) -- (9.5, -0.85);
      \node at (7.5, - 1.25) {$ t^{m + 1}_{f, n}$};
      \node at (5.5, - 1.25) {$t^{m - 1}_{f, n}$};
      \node at (9.5, - 1.25) {$t^{m + 3}_{f, n}$};
      \draw[dashed] (5.5,-0.75) -- (5.5,2);
      \draw[dashed] (7.5,-0.75) -- (7.5,0);
      \draw[dashed] (9.5,-0.75) -- (9.5,2);
      \draw (6.5, -0.65) -- (6.5, -0.85);
      \draw (8.5, -0.65) -- (8.5, -0.85);
      \node at (6.5, - 1.25) {$ t^{m}_{f, n}$};
      \draw[->] (8.5, -2.1) -- (8.5, -1.65);
      \node at (8.5, -2.35) {refine};
      \node at (8.5, - 1.25) {$ t^{m + 2}_{f, n}$};
      \draw[->] (6.5, -2.1) -- (6.5, -1.65);
      \node at (6.5, -2.35) {refine};
    \end{tikzpicture}
    \caption{An example of preserving the local patch structure
      during the marking procedure: if the time step is refined, the other
      time step  belonging to the same patch will also be refined.}
    \label{local_patches}
  \end{figure}
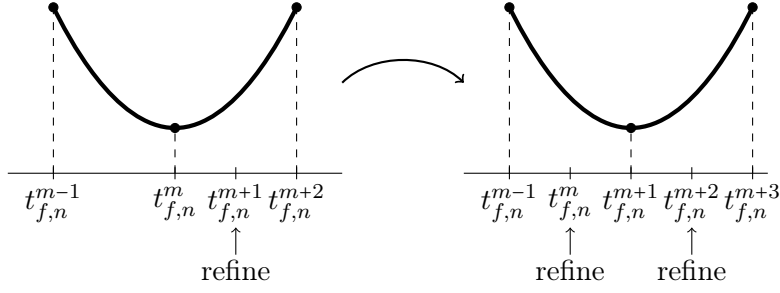

\item Perform time refining.
\item Adjust the macro time-step structure - in case within one macro time-step there exist a fluid and a solid micro time-step that coincide, split the macro time-step into two macro time-steps at this point (see Figure~\ref{splitting}).

  \begin{figure}[t]
    \centering
    \includegraphics[width=\textwidth]{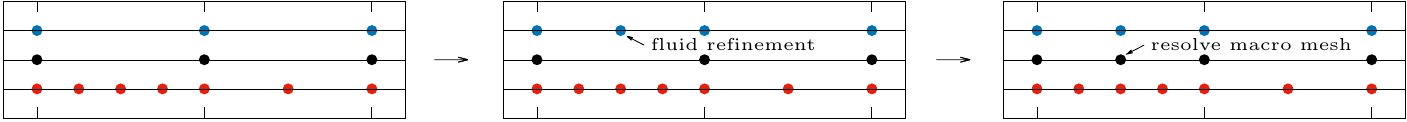}
    \caption{An example of a splitting mechanism of macro
      time-steps. On the left, we show the mesh before refinement:
      middle (in black) the macro nodes, top (in blue) the fluid nodes
      and bottom (in red) the solid nodes with subcycling. In
      the center sketch, we refine the first macro interval once within
      the fluid domain. Since one node is shared between fluid and
      solid, we refine the macro mesh to resolve subcycling. This
      final configuration is shown on the right.}
    \label{splitting}
  \end{figure}

\end{enumerate}


\section{Numerical results}
\label{numerical_results}

\subsection{Fluid subdomain functional}
\label{fluid_functional}

For the first example, we chose to test the derived error estimator on a goal functional concentrated in the fluid subproblem

\begin{equation*}
  J_f(\vec{U})\coloneqq  \int_{0}^T \nu\left(\mathbbm{1}_{\widetilde{\Omega}_f}(\vec{x})\nabla v_f, \nabla v_f\right)_f \diff t, \quad J_s(\vec{U}) \coloneqq 0
\end{equation*}
where $\widetilde{\Omega}_f = (2, 4) \times (0, 1)$ is the right half of the fluid subdomain. For this example, we also took the right hand side concentrated in the fluid subdomain, presented in Configuration~\ref{configuration_1}. As the time interval, we choose $I = [0, 1]$. Then we have
$$(J_f)'_{\vec{U}}(\boldsymbol{\Xi}_f) =  \int_{0}^T 2\nu\left(\mathbbm{1}_{\widetilde{\Omega}_f}(\vec{x})\nabla v_f, \nabla \eta_f\right)_f \diff t.$$
Since the functional is nonlinear, we use a 2-point Gaussian quadrature for integration in time. With~(\ref{midpoints}), the quadrature points read as
\[
g_{f, n}^{m, 1} \coloneqq  \bar t_{f,n}^m + \frac{t_{f, n}^m - t_{f, n}^{m - 1}}{2 \sqrt{3}},\quad
g_{f, n}^{m, 2} \coloneqq  \bar t_{f,n}^m - \frac{t_{f, n}^m - t_{f, n}^{m - 1}}{2 \sqrt{3}}.
\]
With that at hand, we can formulate the discretization of the functional
\begin{equation}\label{functional:fluid}
  \begin{aligned}
  (J_f)'_{\vec{U}}&(\Xi_{f,k}) = \sum_{n = 1}^N\sum_{m = 1}^{M_n}\sum_{q=1}^2
    \nu\left(\mathbbm{1}_{\widetilde{\Omega}_f}(\vec{x}) j_{f,n}^m\nabla v_{f, k}(g_{f, n}^{m, q})j_{f,n}^m\nabla\eta_{f, k}(g_{f, n}^{m, q})\right)_f \\
    &= \sum_{n = 1}^N\sum_{q=1}^2
    \bigg\{\nu(-g_{f, n}^{1, q} + t_{f, n}^1)\left(\mathbbm{1}_{\widetilde{\Omega}_f}(\vec{x}) j_{f,n}^1\nabla v_{f, k}(g_{f, n}^{1, q})\nabla\eta_{f, k}(t_{f, n}^0)\right)_f \\
    &  + \sum_{m = 1}^{M_n - 1} \Big\{\nu(-g_{f, n}^{m + 1, q} + t_{f, n}^{m + 1})\left(\mathbbm{1}_{\widetilde{\Omega}_f}(\vec{x}) j_{f,n}^{m + 1}\nabla v_{f, k}(g_{f, n}^{m + 1, q})\nabla\eta_{f, k}(t_{f, n}^m)\right)_f \\
    & \qquad  \qquad+ \nu(g_{f, n}^{m, q} - t_{f, n}^{m - 1})\left(\mathbbm{1}_{\widetilde{\Omega}_f}(\vec{x}) j_{f,n}^{m}\nabla v_{f, k}(g_{f, n}^{m, q})\nabla\eta_{f, k}(t_{f, n}^m)\right)_f\Big\} \\
  &+ \nu(g_{f, n}^{M_n, q} - t_{f, n}^{M_n - 1})\left(\mathbbm{1}_{\widetilde{\Omega}_f}(\vec{x}) j_{f,n}^{M_n}\nabla v_{f, k}(g_{f, n}^{M_n, q})\nabla\eta_{f, k}(t_{f, n}^{M_n})\right)_f 
  \bigg\},
  \end{aligned}
\end{equation}
where the nodal interpolation is defined as:
\[
j_{f, n}^m\nabla v_{f, k}(t) \coloneqq \frac{t_{f, n}^m - t}{k_{f, n}^m}
\nabla v_{f, k}(t_{f, n}^{m - 1}) + \frac{t - t_{f, n}^{m - 1}}{k_{f,
    n}^m}\nabla v_{f, k}(t_{f, n}^m)
\]
In Table~\ref{fluid_residuals_uniform_equal} we show results of the a posteriori error estimator on a sequence of uniform time meshes. Here, we considered the case without any micro time-stepping, that is the time-step sizes in both fluid and solid subdomains are uniformly equal. That gives a total number of time-steps in the fluid domain equal to $N$ and $N$ in the solid domain. Table~\ref{fluid_residuals_uniform_equal} consists of partial residuals $\theta_{f,k},\theta_{s,k},\vartheta_{f,k}$ and $\vartheta_{s,k}$, overall estimate $\sigma_k$, extrapolated errors $\widetilde{J}-J(\vec U_k)$ and effectivities $\textnormal{eff}_k$. The values of the goal functional on the three finest meshes were used for extrapolation in time. As a result, we got the reference value $\widetilde{J} = 6.029469 \cdot 10^{-5}$. Except for the coarsest mesh, the estimator is very accurate and the effectivities are almost 1. On finer meshes, values of $\theta_{f,k}$ and $\vartheta_{f, k}$ are very close to each other which is due to the linearity of the coupled problem~\cite{BeckerRannacher2001}. A similar phenomenon happens for $\theta_{s, k}$ and $\vartheta_{s, k}$. The residuals are concentrated in the fluid subdomain, which suggests the usage of smaller time-step sizes in this space domain. 

\begin{table}[t]
  \begin{center}
    \resizebox{\textwidth}{!}{%
      \begin{tabular}{c|cccc|ccc}
        \toprule
        $N$  & $\theta_{f, k}$ & $\theta_{s, k}$ & $\vartheta_{f, k}$ & $\vartheta_{s, k}$ &
        $\sigma_{k}$ &
        $\widetilde{J} - J(\vec{U}_k)$ &
        $\textnormal{eff}_k$ \\       
        \midrule
        50 & $3.62\cdot 10^{-8}$ & $5.01 \cdot 10^{-10}$ & $1.05 \cdot 10^{-7}$ &$5.03 \cdot 10^{-10}$     & $1.42 \cdot 10^{-7}$ & $8.06 \cdot 10^{-8}$& 1.76 \\   
        100 & $9.66 \cdot 10^{-9}$ & $1.37 \cdot 10^{-10}$ & $9.96 \cdot 10^{-9}$ & $1.40 \cdot 10^{-10}$  & $1.99 \cdot 10^{-8}$ & $2.05 \cdot 10^{-8}$ & 0.97 \\  
        200 & $2.48 \cdot 10^{-9}$ & $3.00 \cdot 10^{-11}$ & $2.52 \cdot 10^{-9}$ & $3.02 \cdot 10^{-11}$  & $5.07 \cdot 10^{-9}$ & $5.22 \cdot 10^{-9}$ & 0.97 \\  
        400 & $6.28 \cdot 10^{-10}$ & $9.44 \cdot 10^{-12}$ & $6.33 \cdot 10^{-10}$ & $9.56 \cdot 10^{-12}$& $1.28 \cdot 10^{-9}$ & $1.31 \cdot 10^{-9}$ & 0.98 \\  
        800 & $1.58 \cdot 10^{-10}$ & $2.02 \cdot 10^{-12}$ & $1.58 \cdot 10^{-10}$ & $2.06 \cdot 10^{-12}$& $3.20 \cdot 10^{-10}$ & $3.28 \cdot 10^{-10}$ & 0.98 \\
        \bottomrule
      \end{tabular}
    }
    \caption{Residuals and effectivities for fluid subdomain functional in case of uniform time-stepping in case $M_n, L_n = 1$ for all $n$.}
    \label{fluid_residuals_uniform_equal}
  \end{center}
\end{table}

Table~\ref{fluid_residuals_uniform_refined} collects results for another sequence of uniform time meshes. In this case, each of the macro time-steps in the fluid domain is split into two micro time-steps of the same size. That results in $2N$ time-steps in the fluid domain and $N$ in the solid domain. The performance is still highly satisfactory. The residuals remain mostly concentrated in the fluid subdomain. Additionally, after comparing Tables~\ref{fluid_residuals_uniform_equal} and~\ref{fluid_residuals_uniform_refined}, one can see that corresponding values of $\theta_{f, k}$ and $\vartheta_{f, k}$ are the same (value for $N = 800$ in Table~\ref{fluid_residuals_uniform_equal} and $N = 400$ in Table~\ref{fluid_residuals_uniform_refined}, etc.). Overall, introducing micro time-stepping improves performance and reduces extrapolated error $\widetilde{J} - J(\vec{U}_k)$ more efficiently.    

\begin{table}[t]
  \begin{center}
    \resizebox{\textwidth}{!}{%
      \begin{tabular}{c|cccc|ccc}
        \toprule
        $N$  &$\theta_{f, k}$ &
        $\theta_{s, k}$ &
        $\vartheta_{f, k}$ &
        $\vartheta_{s, k}$ &
        $\sigma_{k}$ &
        $\widetilde{J} - J(\vec{U}_k)$&
        $\textnormal{eff}_k$ \\ 
        \midrule
        50 & $9.66 \cdot 10^{-9}$ & $4.99 \cdot 10^{-10}$ & $9.96 \cdot 10^{-9}$ &$5.01 \cdot 10^{-10}$  &$2.06 \cdot 10^{-8}$ & $2.17 \cdot 10^{-8}$& 0.95\\     
        100 & $2.48 \cdot 10^{-9}$ & $1.37 \cdot 10^{-10}$ & $2.52 \cdot 10^{-9}$ & $1.39 \cdot 10^{-10}$&$5.28 \cdot 10^{-9}$ & $5.45 \cdot 10^{-9}$ & 0.97 \\   
        200 & $6.28 \cdot 10^{-10}$ & $2.99 \cdot 10^{-11}$ & $6.33 \cdot 10^{-10}$ & $3.01 \cdot 10^{-11}$&$1.32 \cdot 10^{-9}$ & $1.43 \cdot 10^{-9}$ &  0.92 \\  
        400 & $1.58 \cdot 10^{-10}$ & $9.44 \cdot 10^{-12}$ & $1.58 \cdot 10^{-10}$ & $9.56 \cdot 10^{-12}$&$3.35 \cdot 10^{-10}$ & $3.58 \cdot 10^{-10}$ & 0.94 \\
        \bottomrule
      \end{tabular}
    }
    \caption{Residuals and effectivities for fluid subdomain functional in case of uniform time-stepping in case $M_n = 2$ and $L_n = 1$ for all $n$.}
    \label{fluid_residuals_uniform_refined}
  \end{center}
\end{table}

\begin{table}[t]
  \begin{center}
    \resizebox{\textwidth}{!}{%
    \begin{tabular}{ccc|cccc|ccc}
      \toprule
      $N$ & $M$ & $L$ & $\theta_{f, k}$ & $\theta_{s, k}$ & $\vartheta_{f, k}$ & $\vartheta_{s, k}$&
      $\sigma_{k}$ & $\widetilde{J} - J(\vec{U}_k)$ &$\textnormal{eff}_k$\\
      \midrule
      50 & 56 & 50 & $3.08 \cdot 10^{-8}$ & $5.01 \cdot 10^{-10}$ & $3.16 \cdot 10^{-8}$ & $5.04 \cdot 10^{-10}$ & $6.34 \cdot 10^{-8}$ & $6.64 \cdot 10^{-8}$ & 0.95\\
      50 & 100 & 50 & $9.66 \cdot 10^{-9}$ & $4.99 \cdot 10^{-10}$ & $9.96 \cdot 10^{-9}$ & $5.01 \cdot 10^{-10}$& $2.06 \cdot 10^{-8}$ & $2.17 \cdot 10^{-8}$ & 0.95\\
      50 & 110 & 50 & $8.21 \cdot 10^{-9}$ & $4.99 \cdot 10^{-10}$ & $8.32 \cdot 10^{-9}$ & $5.02 \cdot 10^{-10}$& $1.75 \cdot 10^{-8}$ & $1.84 \cdot 10^{-8}$ & 0.95\\
      50 & 156 & 50 & $5.08 \cdot 10^{-9}$ & $4.99 \cdot 10^{-10}$ & $5.18 \cdot 10^{-9}$ & $4.97 \cdot 10^{-10}$ & $1.13 \cdot 10^{-8}$ & $1.20 \cdot 10^{-8}$ & 0.94\\
      \bottomrule
    \end{tabular}}
    \caption{Residuals and effectivities for fluid subdomain functional in case of adaptive time-stepping.}
    \label{fluid_residuals_adaptive}
  \end{center}
\end{table}

In Table \ref{fluid_residuals_adaptive} we present findings in the case of adaptive time mesh refinement. We chose an initial configuration of uniform time-stepping without micro time-stepping for $N = 50$ and applied a sequence of adaptive refinements. On every level of refinement, the total number of time-steps is $M + L$. One can see that since the error is concentrated in the fluid domain, only time-steps corresponding to this space domain were refined. Again, effectivity gives very good results. The extrapolated error $\widetilde{J} - J(\vec{U}_k)$ is even more efficiently reduced. 

\subsection{Solid subdomain functional}
\label{solid_functional}

For the sake of symmetry, for the second example, we chose a functional concentrated on the solid subdomain
\begin{equation*}
  J_f(\vec{U}) = 0, \quad J_s(\vec{U}) =  \int_{0}^T \lambda\left(\mathbbm{1}_{\widetilde{\Omega}_s}(\vec{x})\nabla u_s, \nabla u_s\right)_s \diff t,
\end{equation*}
where $\widetilde{\Omega}_s = (2, 4) \times (-1, 0)$ is the right half of the solid subdomain. This time we set the right hand side according to Configuration~\ref{configuration_2}. Again, $\bar{I} = [0, 1]$. The derivative reads as
\[
(J_s)'_{\vec{U}}(\boldsymbol{\Xi}_s) =  \int_{0}^T 2\lambda\left(\mathbbm{1}_{\widetilde{\Omega}_s}(\vec{x})\nabla u_s, \nabla \xi_s\right)_s \diff t,
\]
and allows for a discretization according to~(\ref{functional:fluid}). 
Similarly, Table~\ref{solid_residuals_uniform_equal} gathers results for a sequence of uniform meshes without any micro time-stepping ($N + N$ micro time-steps). The last three solutions are used for extrapolation in time which gives $\widetilde{J} = 3.458826 \cdot 10^{-4}$. Also for this example, the effectivity is very satisfactory. On the finest discretization, the effectivity slightly declines. This might come from the limited accuracy of the reference value. Once more, on finer meshes, fluid residuals $\theta_{f, k}$, $\vartheta_{f, k}$ and solid residuals $\theta_{s, k}$ $\vartheta_{s, k}$ have similar values. This time, the residuals are concentrated in the solid subdomain and, in this case, the discrepancy is a bit bigger.  

\begin{table}[t]
  \begin{center}
    \resizebox{\textwidth}{!}{
      \begin{tabular}{c|cccc|ccc}
        \toprule
        $N$  & $\theta_{f, k}$ & $\theta_{s, k}$ &$\vartheta_{f, k}$ & $\vartheta_{s, k}$&$\sigma_{k}$ &$\widetilde{J} - J(\vec{U}_k)$ &$\textnormal{eff}_k$ \\
        \midrule
        50 & $2.03 \cdot 10^{-10}$ & $2.66 \cdot 10^{-6}$ & $1.93 \cdot 10^{-10}$ &$1.03 \cdot 10^{-5}$  & $1.30 \cdot 10^{-5}$ & $2.49 \cdot 10^{-5}$& 0.52 \\ 
        100 & $4.53 \cdot 10^{-11}$ & $2.59 \cdot 10^{-6}$ & $4.26 \cdot 10^{-11}$ & $2.67 \cdot 10^{-6}$& $5.26 \cdot 10^{-6}$ & $4.77 \cdot 10^{-6}$ & 1.10 \\
        200 & $1.28 \cdot 10^{-11}$ & $5.18 \cdot 10^{-7}$ & $1.26 \cdot 10^{-11}$ & $5.21 \cdot 10^{-7}$& $1.04 \cdot 10^{-6}$ & $9.80 \cdot 10^{-7}$ & 1.06 \\
        400 & $3.30 \cdot 10^{-12}$ & $1.17 \cdot 10^{-7}$ & $3.29 \cdot 10^{-12}$ & $1.17 \cdot 10^{-7}$& $2.34 \cdot 10^{-7}$ & $2.23 \cdot 10^{-7}$ & 1.05 \\
        800 & $8.32 \cdot 10^{-13}$ & $2.82 \cdot 10^{-8}$ & $8.32 \cdot 10^{-13}$ & $2.80 \cdot 10^{-8}$& $5.62 \cdot 10^{-8}$ & $5.07 \cdot 10^{-8}$ & 1.11 \\
        \bottomrule
      \end{tabular}}
    \caption{Residuals and effectivities for solid subdomain functional in case of uniform time-stepping in case $M_n, L_n = 1$ for all $n$.}
    \label{solid_residuals_uniform_equal}
  \end{center}
\end{table}

\begin{figure}[t]
  \begin{center}
    \includegraphics[width=\textwidth]{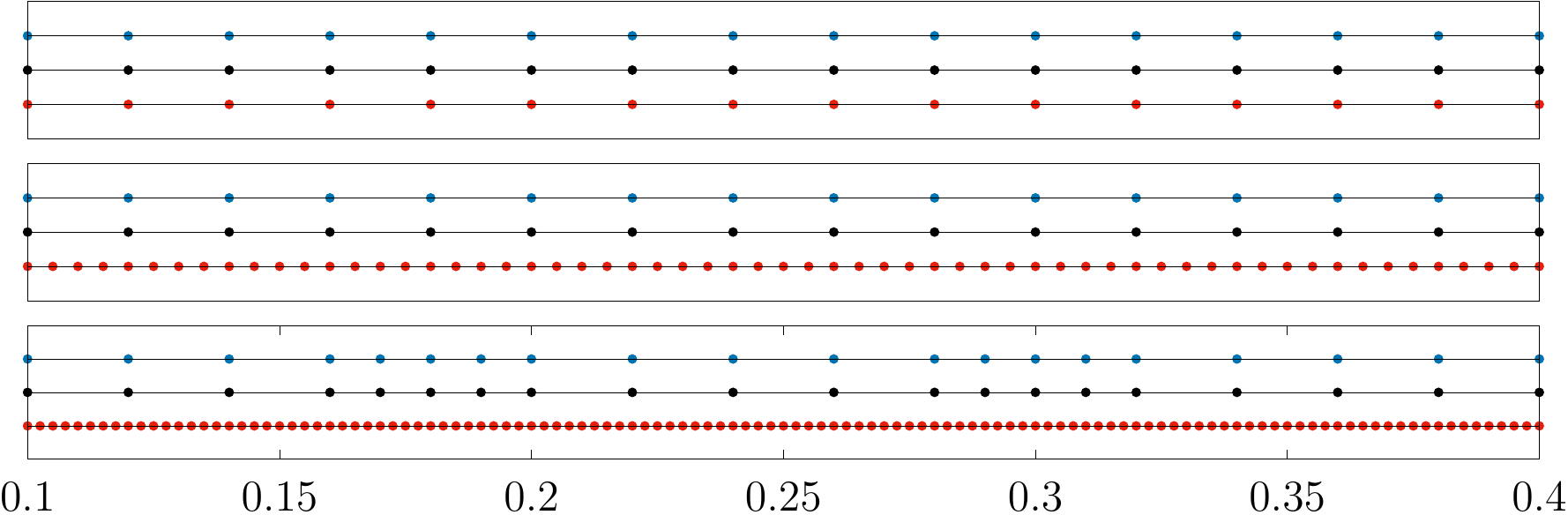}
  \end{center}
  \caption{Adaptive meshes for the solid functional. Top: uniform
    initial mesh; middle: 2 steps of adaptive refinement; bottom: 4
    steps. Each plot shows the macro mesh (middle), the fluid mesh
    (top, in blue) and the solid mesh (bottom, in red). }
  \label{fig:refine}
\end{figure}

In Table~\ref{solid_residuals_uniform_refined} we display outcomes for a sequence of uniform meshes where each of the macro time-steps in the solid subdomain is split into two micro time-steps. That gives $N + 2N$ time-steps. Introducing micro time-stepping does not have a negative impact on the effectivity and significantly saves computational effort.
Corresponding values of $\theta_{s,k}$ and $\vartheta_{s, k}$ in Tables \ref{solid_residuals_uniform_equal} and 
\ref{solid_residuals_uniform_refined} are almost the same. Residuals remain mostly concentrated in the solid subdomain. 

\begin{table}[t]
  \begin{center}
    \resizebox{\textwidth}{!}{%
      \begin{tabular}{c|cccc|ccc}
        \toprule
        $N$  &$\theta_{f, k}$ &$\theta_{s, k}$ & $\vartheta_{f, k}$ &$\vartheta_{s, k}$&$\sigma_{k}$ & $\widetilde{J} - J(\vec{U}_k)$ &  $\textnormal{eff}_k$ \\
        \midrule
        50 & $4.13 \cdot 10^{-10}$ & $2.61 \cdot 10^{-6}$ & $1.91 \cdot 10^{-9}$ &$2.68 \cdot 10^{-6}$ & $5.29 \cdot 10^{-6}$ & $4.68 \cdot 10^{-6}$& 1.13\\  
        100 & $8.69 \cdot 10^{-11}$ & $5.20 \cdot 10^{-7}$ & $-3.72 \cdot 10^{-11}$ & $5.23 \cdot 10^{-7}$& $1.04 \cdot 10^{-6}$ & $9.54 \cdot 10^{-7}$ & 1.09 \\
        200 & $1.80 \cdot 10^{-11}$ & $1.17 \cdot 10^{-7}$ & $1.40 \cdot 10^{-12}$ & $1.17 \cdot 10^{-7}$ & $2.34 \cdot 10^{-7}$ & $2.16 \cdot 10^{-7}$ & 1.08 \\
        400 & $3.94 \cdot 10^{-12}$ & $2.82 \cdot 10^{-8}$ & $1.87 \cdot 10^{-12}$ & $2.80 \cdot 10^{-8}$ & $5.62 \cdot 10^{-8}$ & $4.90 \cdot 10^{-8}$ & 1.15 \\
        \bottomrule
      \end{tabular}}
    \caption{Residuals and effectivities for solid subdomain functional in case of uniform time-stepping in case $M_n = 1$ and $L_n = 2$ for all $n$.}
    \label{solid_residuals_uniform_refined}
  \end{center}
\end{table}

Following the fluid example, in Table~\ref{solid_residuals_adaptive} we show calculation results in the case of adaptive time mesh refinement. Here as well we took the uniform time-stepping without micro time-stepping for $N = 50$ as the initial configuration and the total number of time-steps is $M + L$. Except for the last entry, only the time-steps corresponding to the solid domain were refined. On the finest mesh, the effectivity deteriorates. However, adaptive time-stepping is still the most effective in reducing the extrapolated error $\widetilde{J} - J(\vec{U}_k)$.

\begin{table}[t]
  \begin{center}
    \resizebox{\textwidth}{!}{%
      \begin{tabular}{ccc|cccc|ccc}
        \toprule
        $N$ & $M$ & $L$ & $\theta_{f, k}$ &$\theta_{s, k}$ &$\vartheta_{f, k}$ &$\vartheta_{s, k}$&
        $\sigma_{k}$ &$\widetilde{J} - J(\vec{U}_k)$ & $\textnormal{eff}_k$\\
        \midrule
        50 & 50 & 88 & $3.77 \cdot 10^{-10}$ & $6.57 \cdot 10^{-6}$ & $6.72 \cdot 10^{-8}$ & $6.91 \cdot 10^{-6}$  & $1.35 \cdot 10^{-5}$ & $1.06 \cdot 10^{-5}$ & 1.28\\
        50 & 50 & 166 & $5.17 \cdot 10^{-10}$ & $1.35 \cdot 10^{-6}$ & $7.16 \cdot 10^{-8}$ & $1.38 \cdot 10^{-6}$ & $2.80 \cdot 10^{-6}$ & $2.52 \cdot 10^{-6}$ & 1.11\\
        50 & 50 & 286 & $5.80 \cdot 10^{-10}$ & $4.54 \cdot 10^{-7}$ & $4.16 \cdot 10^{-8}$ & $4.56 \cdot 10^{-7}$ & $9.52 \cdot 10^{-7}$ & $7.34 \cdot 10^{-7}$ & 1.30\\
        54 & 54 & 400 & $5.70 \cdot 10^{-10}$ & $1.19 \cdot 10^{-7}$ & $4.12 \cdot 10^{-8}$ & $1.19 \cdot 10^{-7}$ & $2.81 \cdot 10^{-7}$ & $1.10 \cdot 10^{-7}$ & 2.55\\
        \bottomrule
      \end{tabular}}
    \caption{Residuals and effectivities for solid subdomain functional in case of adaptive time-stepping.}
    \label{solid_residuals_adaptive}
  \end{center}
\end{table}

Finally, we show in Figure~\ref{fig:refine} a sequence of adaptive
meshes that result from this adaptive refinement strategy. In the top
row, we show the initial mesh with 50 macros steps and no further
splitting in fluid and solid. For a better presentation, we only show a
small subset of the temporal interval $[0.1,0.4]$. In the middle plot,
we show the mesh after 2 steps of adaptive refinement and in the
bottom line after 4 steps of adaptive refinement. Each plot shows the
macro mesh, the fluid mesh (above) and the solid mesh (below). As
expected, this example leads to a sub-cycling within the solid
domain. For a finer approximation, the fluid problem also requires
some local refinement. Whenever possible we avoid excessive subcycling
by refining the macro mesh as described in
Section~\ref{adaptivity}.

%

\section{Conclusion}
In this paper, we have developed a multirate scheme and a temporal error estimate for a coupled problem that is inspired by fluid-structure interactions. The two subproblems, the heat equation and the wave equation, feature different temporal dynamics such that balanced approximation properties and stability demands ask for different step sizes.  

We introduced a monolithic variational Galerkin formulation for the coupled problem and then used a partitioned framework for solving the algebraic systems. Having different time-step sizes for each of the subproblems couples multiple states in each time-step, which would require an enormous computational effort. To solve this, we discussed two different decoupling methods: first, a simple relaxation scheme that alternates between fluid and solid problem and second, similar to the  shooting method, where we defined a root-finding problem on the interface and used matrix-free Newton-Krylov method for quickly approximating the zero.  Both of the methods were able to successfully decouple our specific example and showed good robustness concerning different subcycling of the multirate scheme in fluid- or solid-domain. However, the convergence of the shooting method was faster and it required fewer evaluations of the variational formulation. 

As the next step, we introduced a goal-oriented error estimate based on the dual weighted residual method to estimate errors with regard to functional evaluations. The monolithic space-time Galerkin formulation allowed to split the residual errors into contributions from the fluid and solid problems. Several numerical results for two different goal functionals show very good effectivity of the error estimate. Finally, we established the localization of the error estimator. That let us derive an adaptive refinement scheme for choosing optimal distinct time meshes for each problem.

In future work, it remains to extend the methodology to nonlinear problems, in particular, to fully coupled fluid-structure interactions. 


\section{Acknowledgements}
Both authors acknowledge support by the Deutsche Forschungsgemeinschaft (DFG, German Research Foundation) - 314838170, GRK 2297 MathCoRe. TR further acknowledge supported by the Federal Ministry of Education and Research of Germany (project number 05M16NMA).

\bibliographystyle{ieeetr}
\bibliography{lit}

\end{document}